\documentclass[12pt,reqno]{amsart}
\advance \topmargin by -\headheight
\advance \topmargin by -\headsep
\evensidemargin \oddsidemargin
\usepackage{amsmath, amsthm, amssymb, amsfonts}
\usepackage{mathtools,thmtools}
\usepackage{hyperref}
\usepackage[capitalise]{cleveref}
\usepackage{enumerate}
\usepackage{url}
\usepackage{dsfont}
\usepackage{mathrsfs}
\DeclareMathOperator{\tr}{tr}

\newtheorem{theorem}{Theorem}[section]
 
\newtheorem{lemma}[theorem]{Lemma}
\newtheorem{corollary}[theorem]{Corollary}

\usepackage{bm}

\usepackage{xcolor}

\usepackage{fancyvrb}

\newtheorem{prop}[theorem]{Proposition}

\theoremstyle{definition}
\newtheorem{defn}[theorem]{Definition}
\newtheorem{example}[theorem]{Example}

\theoremstyle{remark}
\newtheorem{remark}[theorem]{Remark}

\numberwithin{equation}{section}

\renewcommand{\geq}{\geqslant}
\renewcommand{\leq}{\leqslant}

\def\eps{\varepsilon}

\newcommand*{\EE}{
	\mathop{
		\mathchoice{\vcenter{\hbox{\larger[4]$\mathbb{E}$}}}
		{\kern0pt\mathbb{E}}
		{\kern0pt\mathbb{E}}
		{\kern0pt\mathbb{E}}
	}\displaylimits
}

\def \C {\mathbb{C}}
\def \N {\mathbb{N}}
\def \Q {\mathbb{Q}}
\def \R {\mathbb{R}}

\def \Z {\mathbb{Z}}

\DeclareMathOperator{\supp}{supp}

\newcommand\1{\mathds{1}}

\DeclarePairedDelimiter{\abs}{\lvert}{\rvert}
\DeclarePairedDelimiter{\n}{\lVert}{\rVert}

\title[Additive energy, Freiman's lemma, and edge-isoperimetry]{Additive energy, Freiman's lemma, and an edge-isoperimetric problem on the boolean cube}

\author{Akshat Mudgal}
\email{Akshat.Mudgal@warwick.ac.uk}
\address{Mathematics Institute, Zeeman Building, University of Warwick, Coventry CV4 7AL, UK}

\author{Carl Schildkraut}
\email{carlsch@stanford.edu}
\address{Department of Mathematics, Stanford University, Stanford, CA 94305, USA}

\author{Dmitrii Zakharov}
\email{zakhdm@mit.edu}
\address{Department of Mathematics, Massachusetts Institute of Technology, Cambridge, MA 02139, USA}


\begin{document}
	
	\begin{abstract}
		Let $d \in \mathbb{N}$ and let $A \subseteq \mathbb{R}^d$ be a finite, non-empty set not contained in any translate of a hyperplane. 
		We show that there exists a subset $S \subseteq A$ such that $d \ll |S| \ll d$ and
		\[E(A,S) := |\{(a_1,a_2,s_1,s_2)\in A\times A\times S\times S:a_1+s_1=a_2+s_2\}|\ll\lvert A\rvert\lvert S\rvert.\]
		This can be interpreted as a robust, analytic generalisation of the classical Freiman's lemma, and recovers Freiman's lemma up to a constant multiplicative factor. In the process of proving this result, we also show that for any integer $1 \leq m \leq d$ and any set $A \subset \mathbb{R}^d$ with $|A| \leq \binom{d}{m}$, one has
		\[ E(A, B) \leq 4(m + 1/m) |A||B|,\]
		where $B$ is the canonical basis of $\mathbb{Z}^d$. This gives new estimates for multiplicative energies between finite sets of primes and arbitrary subsets of positive integers, improving a previous result of Shkredov. Using related ideas, we also prove new bounds for an edge-isoperimetric problem for certain graphs on the boolean cube, improving an earlier result of Kirshner and Samorodnitsky. We give two different proofs of the first result described in this abstract and three different proofs of the second one.
	\end{abstract}
	
	\maketitle
	
	\section{Introduction} \label{introductionsection}

	An important theme in additive combinatorics is the interplay between notions of dimension of a set and its additive structure. For example, Freiman's theorem states that any finite set $A \subseteq \mathbb{Z}$ such that $\abs{A+A}\leq K\abs{A}$, where 
	\[ A+A = \{a_1 + a_2 : a_1, a_2 \in A\},\]
	must be contained in a small-dimensional \emph{generalised arithmetic progression} $P$ of the form
	\[ P = \{ u_0 +  u_1 \ell_1 + \dots + u_r \ell_r :  1 \leq \ell_i \leq L_i \ \ (1 \leq i \leq r)  \}\]
	for some suitably chosen $u_0, u_1, \dots, u_r \in \mathbb{Z}$ and $L_1, \dots, L_r \in \mathbb{N}$ satisfying 
	\begin{equation}  \label{conc1}
		L_1 \cdots L_r \ll_K |A| \ \ \text{and} \ \ r \ll_K 1.
	\end{equation} 
	Obtaining variations of this type of result with a strong quantitative dependence on $K$ forms the content of the infamous polynomial Freiman--Ruzsa conjecture over $\mathbb{Z}$, see \cite{Ra2025,Sa2013}.
	
	From the perspective of applications, one is often just interested in  obtaining good bounds for the second inequality in \eqref{conc1}; see, for instance, the weak polynomial Freiman--Ruzsa conjecture over $\mathbb{Z}$ and its applications to sum-product-type estimates over $\mathbb{R}$  \cite{GGMT2025, Mu2024}. This is roughly equivalent to studying the following type of question: given a set $A \subseteq \mathbb{R}^d$, what notions of high-dimensionality can we apply on $A$ to deduce strong lower bounds for $|A+A|/|A|$? The starting result here is the well-known Freiman's lemma, which states that  if $A \subseteq \mathbb{R}^d$ satisfies $\dim(A) = d$,  then 
	\begin{equation} \label{frlem}
		|A+A| \geq (d+1)|A| - d(d+1)/2 \geq (d+1)|A|/2,
	\end{equation}
	where $\dim(A)$ denotes the dimension of the affine subspace generated by $A$. The first inequality here is sharp, as can be seen by setting 
	\[ A = \{ e_1, 2e_1, \dots, Ne_1\} + \{0, e_2, e_3, \dots, e_d\} , \]
	where $\{e_1, \dots, e_d\}$ is the canonical basis of $\mathbb{Z}^d$. One can obtain stronger lower bounds of the shape $|A+A| \geq C^d|A|$ for some absolute constant $C>1$ by assuming stricter ``higher-dimensional" conditions on $A$; see, for example, the Freiman--Bilu theorem \cite{Bilu1999} or lower bounds for sumsets of sets containing hypercubes  \cite{GT2006, LMPS2025}. Despite its simplicity, \eqref{frlem} has had various applications to Freiman-type inverse theorems and sum-product estimates. 
	
	More recently, various works \cite{JM2023, CDM2024} have explored robust versions of \eqref{frlem}, wherein, given a finite set $A \subseteq \mathbb{R}^d$ with $\dim(A) =d$, one is interested in finding a set $S \subseteq A$ such that 
	\begin{equation} \label{rbfr}
		|S| \ll_d 1 \ \ \text{and} \ \ |A+S| \gg d|A|. 
	\end{equation}
	Such results have  found applications towards a conjecture of Alon \cite{Al2007} on independence numbers of Cayley graphs over $\mathbb{Z}/N\mathbb{Z}$  as well as a question of Green--Morris \cite{GM2016} on counting the number of subsets of $[N] = \{1,2\dots,N\}$ with a given size and doubling, see \cite{CDM2024, CDM2026}. Jing--Mudgal \cite{JM2023} proved that for any finite set $A \subseteq \mathbb{R}^d$ with $\dim(A) = d$, one can find $S \subseteq A$ with 
	\begin{equation} \label{bns2}
		|S| \ll_d 1 \ \ \text{and}  \ \ |A+S| \geq (d+1)|A| -5(d+1)^3. 
	\end{equation}
	The dependence of $|S|$ on $d$ here is quite poor.
    By employing quite different methods, Campos--Dahia--Marciano \cite{CDM2024} proved that for any $\gamma >0$, whenever $d > \gamma^{-3}$, one can find $S \subseteq A$ such that 
	\begin{equation} \label{bns3}
		|S| \ll \gamma^{-1} d \ \ \text{and} \ \ |A+S| \geq (1 - 5\gamma)(d+1)|A|/2.  
	\end{equation}
	Thus, while the above lower bound is weaker than its analogue described in \eqref{bns2}, the upper bound on $|S|$ is quite nice. In a very recent sequel paper, Campos--Dahia--Marciano \cite{CDM2026} removed the factor of $1/2$ from the lower bound in the right hand side of \eqref{bns3} at the cost of worsening the dependence of $\abs{S}$ on $d$ as well as requiring $\abs{A} \gg d/\gamma$.

	In this paper, using very different techniques as compared to \cite{CDM2024, CDM2026, JM2023} we provide a significant qualitative generalisation of such results which is more analytic in nature.  In particular, let us define, for any non-empty finite sets $A,B \subseteq \mathbb{R}^d$, the \emph{additive energy}
	\[ E(A,B) = \sum_{a_1, a_2 \in A} \sum_{b_1, b_2 \in B} \1_{a_1 + b_1 = a_2 + b_2} .  \]
	Moreover, we call a non-empty set $S = \{s_1, \dots, s_r\} \subset \R^d$ \emph{affine-independent} if, for any $c_1, \dots, c_r \in \R$ with $c_1 + \dots + c_r =0$, we have $c_1 s_1 + \dots + c_r s_r =0$ if and only if $c_1 = \dots = c_r = 0$.
	With this in hand, our main result is as follows.

	\begin{theorem}\label{thm:main}
		Let $d \in \mathbb{N}$ and let $A\subset\R^d$ be a finite, non-empty set with $\dim(A) = d$. Then there exists some affine-independent $S\subset A$ such that $d \ll |S| \ll d$ and
		\[E(A,S)\ll \abs{A}\abs{S}.\]
	\end{theorem}

	\cref{thm:main} immediately delivers estimates of the shape \eqref{rbfr}. Indeed, a simple application of Cauchy's inequality implies that
	\[ |A+S| \geq |A|^2|S|^2/E(A,S) \gg |A||S| \gg d|A|.  \]
	Hence, Theorem \ref{thm:main} suggests a new heuristic as to why finite sets $A \subset \mathbb{R}^d$ with $\dim(A) =d$ satisfy $|A+A|/|A| \gg d$: one can always find an affine-independent set $S \subseteq A$ satisfying $|S| \gg d$ with which $A$ has very little additive interaction. More generally, Theorem \ref{thm:main} dispenses the following corollary.
	
	\begin{corollary} \label{cor56}
		Let  $\eta\in (0,1]$, let $d \in \mathbb{N}$ and let $A\subset\R^d$ be a finite, non-empty set with $\dim(A) = d$. Then there exists some $S\subset A$ with $d \ll |S| \ll d$  such that, for any set $B \subseteq A$ with $|B| \geq \eta|A|$, one has
		\[ |B+S| \gg \eta^2 d |A|. \]
	\end{corollary}
	
	This highlights the strength of upper bounds on energies as compared to lower bounds of sumsets; the former very easily extend to dense subsets while the latter do not offer this flexibility. Consequently, our proof of Theorem \ref{thm:main} requires some new ideas and in fact, the techniques involved therein differ significantly from the ideas used in \cite{JM2023, CDM2024, CDM2026}. Indeed, while the latter works rely crucially on combinatorial-geometric arguments to obtain translates of the set $A$ which are sufficiently disjoint, one of the ideas we use to obtain upper bounds on the additive energy $E(A,S)$ is to employ analytic methods to pass to studying higher-order energies.
	
	Thus, given finite, non-empty subsets $X,S$ of some abelian group and $f\colon X\to[0,\infty)$, we define the \emph{weighted additive energy}
	\begin{equation}  \label{wten}
		E_f(X,S):=\sum_{x_1,x_2\in X}f(x_1)f(x_2)\sum_{s_1,s_2\in S}\1_{x_1+s_1=x_2+s_2}.
	\end{equation}
	We are interested in upper bounds for $E_f(X,S)$ and $E(X,S)$ when $S = \{e_1, \dots, e_d\}$ and $X \subseteq \mathbb{R}^d$ is some finite, non-empty set. When $|X| \geq 2^d$, one can not do significantly better than the trivial bound $E(X,S) \leq |X||S|^2$. Indeed setting $X = \{0,1\}^d$, observe that for any $1 \leq i < j \leq d$, there are at least $2^{d-2}$ elements $x \in X$ such that $x + e_i - e_j \in X$, and so,
	\[ E(X,S) \geq |S||X| +  |S|(|S|-1) |X|/4 \gg |X||S|^2. \]
	On the other hand, when $|X| \ll d^{O(1)}$, one can not find such good lower bounds. For instance, given $d, m \in \mathbb{N}$ with $m \leq d/2$, one can set
	\[ X = \{ (x_1, \dots, x_d) \in \{0,1\}^d : x_1 + \dots + x_d = m  \}. \]
	For any $x \in X$, there are at least $m(d-m)$ choices of $1 \leq i, j \leq d$ such that $x - e_i + e_j \in X$, whence,
	\begin{equation} \label{str12}
		E(X,S) \geq m(d-m) |X|. 
	\end{equation}
	We are able to prove a corresponding upper bound in this setting.
	
	\begin{theorem} \label{asymenergy}
		Let $d,m \in \mathbb{N}$ with $m \leq d$, let $S = \{e_1, \dots, e_d\}$ and let 
		$f\colon \R^d\to[0,\infty)$ be supported on a finite, non-empty set $X \subseteq \mathbb{R}^d$ and satisfy
		\begin{equation}\label{ell12}
			\frac{\left(\sum_{x\in X}f(x)\right)^2}{\sum_{x\in X}f(x)^2}\leq\binom dm.
		\end{equation}
		Then we have
		\begin{equation} \label{srs}
			E_f(X,S)\leq 4d(m + 1/m)  \sum_{x \in X} f(x)^2 .
		\end{equation}
	\end{theorem}

	Setting $f(x) = \1_X$ for some $X\subset\R^d$ with $|X| \leq \binom{d}{m}$, \cref{asymenergy} implies
	\begin{equation} \label{4jl}
		E(X,S) \ll md |X|,
	\end{equation}
	matching the lower bound in \eqref{str12} up to a factor of $O(1)$.

	It turns out that such estimates are very closely related to an edge-isoperimetric problem for graphs on the boolean cube.
	For any set $X\subseteq\{0,1\}^d$, write $T(X)$ for the number of ordered pairs $(x_1,x_2)\in X\times X$ for which the difference $x_1-x_2$ contains exactly two non-zero entries. Thus, $T(X)$ counts the number of pairs of elements of $X$ at Hamming distance~$2$. The quantity $T(X)$ has been previously studied in \cite{ACKS2005, KS2021}, and the best known bound here is due to Kirshner--Samorodnitsky \cite[Corollary~7]{KS2021}, who proved strong hypercontractivity-type inequalities and combined them with Fourier-analytic arguments to get
	\begin{equation}   \label{ksresult}
		T(X) \leq e^2  m(d-m)|X| \ \ \text{whenever} \ X \subseteq \{0,1\}^d \ \text{satisfies} \ |X| \leq \binom{d}{m} 
	\end{equation}
	for some integer $1 \leq m \leq d/2$. The trivial upper bound on $T(X)$ is $\binom d2\abs{X}$ (see \eqref{atrivbd}), and so \eqref{ksresult} is only nontrivial when $m\leq 0.073d$. 
	
	Using a more direct graph-theoretic approach, we can prove the following upper bound on $T(X)$, which improves upon \eqref{ksresult} whenever either bound is nontrivial.

	\begin{theorem}\label{hammingcor}
		Let $m,d$ be integers such that $1 \leq m \leq d/2$ and let $X\subset\{0,1\}^d$ satisfy $\abs{X}\leq\binom dm$. Then
		\begin{equation}\label{imp3}
			T(X)\leq 4md\abs{X}.
		\end{equation}
	\end{theorem}
	
	We highlight the connection between \cref{asymenergy,hammingcor} in Remark \ref{alprem}. For now, we simply note that a bound of the form \eqref{4jl} for $X\subseteq\{0,1\}^d$ follows directly by applying \cref{hammingcor}, since points $x_1,x_2\in X$ satisfying $x_1+e_i=x_2+e_j$ for distinct $1 \leq i,j \leq d$ are at Hamming distance $2$. Our combinatorial proof of \cref{hammingcor} can also be used to prove \eqref{4jl} in general with the same leading constant $4$; see \cref{rmk:gen}.

	Another interesting aspect of the proofs of Theorems \ref{asymenergy} and \ref{hammingcor} is that they do not invoke any hypercontractivity-type result or Fourier analysis. 
	That said, we will provide an alternative approach to proving a weaker version of Theorem \ref{asymenergy} using results of Kirshner--Samorodnitsky \cite{KS2021} on hypercontractivity for the boolean cube. In particular, we will use the latter to prove a variant of \cref{asymenergy} in the group $(\Z/2\Z)^d$, see \cref{lem:ks-gen}.
	This variant allows us to deduce \cref{prop:weighted}, an analogue of \cref{asymenergy} which gives a weaker version of conclusion \eqref{srs} on the stricter assumption that $\abs{X}\leq\binom dm$. We remark that either \cref{asymenergy} or \cref{prop:weighted} is enough to prove \cref{thm:main}.

	Returning to \eqref{4jl}, we can use this along with a standard fibring argument and the theory of Freiman isomorphisms to deduce the following.
	
	\begin{corollary} \label{primeresult}
		Let $d,m \in \mathbb{N}$ such that $m \leq d$, let $P \subseteq \mathbb{N}$ be a set of $d$ prime numbers and let $Q \subseteq \mathbb{N}$ be a finite set such that $|Q| \leq \binom dm$. Then 
		\[ M(P,Q) := \sum_{q_1, q_2 \in Q} \sum_{p_1, p_2 \in P} \1_{p_1 q_1 = p_2 q_2} \ll m d |Q|. \]
	\end{corollary}
	
	Corollary \ref{primeresult} improves upon an earlier result of Shkredov \cite{Shk2022} which states that for any set $P$ of $d$ prime numbers and any finite, non-empty set $Q \subset \mathbb{N}$, one has
	\begin{equation} \label{shkredovres1}
		M(P,Q) \ll d |Q| \log |Q|.  
	\end{equation}
	The above bound is only non-trivial when $|Q| = 2^{o(d)}$.
	For such choices of $d$, let $1 \leq m \leq d/2$ be an integer such that ${d \choose m-1} \leq |Q| \leq {d \choose m}$. Then $\log |Q| \gg (m-1) \log (d/m)$, and so, Corollary \ref{primeresult} improves \eqref{shkredovres1} by a factor of $\log(d/m)$.

	\begin{remark} \label{shk4}
		We note that Shkredov \cite{Shk2022} studied the setting of Corollary \ref{primeresult} in order to prove the  following result. Let $A \subset \mathbb{N}$ be a finite, non-empty set such that $|A \cdot A|  \ll |A|$, where  $A \cdot A =\{a \cdot a' : a, a' \in A\}$ denotes the product set of $A$. Then the smallest positive integer $k \notin A$ satisfies
		\begin{equation}  \label{shkredoveq2} 
			k \ll \log |A| \log \log |A|. 
		\end{equation}
		This can be seen as an integer analogue of Vinogradov's least quadratic nonresidue problem which, given prime $p$, concerns the size $n_p$ of the least $k \in \mathbb{N}$ such that $k \ ({\rm mod} \ p) \notin R$, where $R$ is the set of quadratic residues in $\mathbb{Z}/p\mathbb{Z}$. Note that $|R \cdot R| = |R|$. It is expected that $n_p \ll_{\varepsilon} p^{\varepsilon}$ for every $\varepsilon>0$, but the best unconditional result towards this is of the form $n_p \ll_{\varepsilon} p^{1/(4 e^{1/2}) + \varepsilon}$ for every $\varepsilon >0$.\footnote{This conjecture has now seemingly been resolved by OpenAI in their very recent work on zero-free regions of primitive Dirichlet $L$-functions.} The latter is a classical result of Burgess \cite{Bur1957}.  \

		We  observe that one can use Freiman's lemma directly to prove that 
		\[ k \ll 1. \]
		Indeed, let $\{1,2,\dots, k-1\} \subseteq A$ and let $P$ be the set of primes in $\{1,2,\dots, k-1\}$. One may now use a standard fibring argument along with \eqref{frlem} to deduce that
		\[ |P||A|/2 \leq |A \cdot A| \ll |A|\]
		whence $|P| \ll 1$. Since $|P| \gg k / \log k$, we must have $k \ll 1$.
	\end{remark}

    \subsection*{Proof Ideas} We briefly discuss some of the ideas employed in our paper. Our starting point for proving Theorem \ref{thm:main} is by proving a decomposition lemma for arbitrary $d$-dimensional sets $A \subset \R^d$, see Lemma \ref{lem:decomp}. In particular, using iterative applications of the weak Polynomial Freiman--Ruzsa theorem due to Gowers--Green--Manners--Tao and a strong quantitative version of the Balog--Szemer\'{e}di--Gowers theorem due to Reiher--Schoen, we are able to show that any set $d$-dimensional set $A \subset \R^d$ can be partitioned into a ``structured piece'' $A_{\rm str}$ and an ``unstructured piece'' $A_{\rm unstr}$, where $A_{\rm str}$ is contained in few translates of a sufficiently low-dimensional subspace $V \subset \R^d$ while $A_{\rm unstr}$ either is quite small, or has a small additive energy and is not concentrated on any translate of subspaces approximating $V$. Using variants of Theorem \ref{asymenergy}, we are able to control $E(A_{\rm str}, S)$ for \emph{any} set $S \subset \R^d$ such that the set $\pi(S)$ is affine-independent and $|\pi(S)|=|S|$, where $\pi\colon \R^d \to \R^d/V$ is the natural projection map. Our goal now is to find a suitable affine-independent set $S \subset A$ such that $E(A_{\rm unstr}, S)$ is also well-controlled. We do this via employing the probabilistic method.   More specifically, we use a novel, iterative sampling procedure that takes into account that our random set $S$ must have $\pi(S)$ be affine-independent. The aforementioned properties of $A_{\rm unstr}$ then imply that $E(A_{\rm unstr}, S)$ is suitably bounded above.

    As for Theorem \ref{asymenergy}, we provide three different approaches to proving such estimates, the first being a graph-theoretic approach, the second relying on a technical analysis of higher-order additive energies and the final relying on sharp hypercontractivity-type results. Here, we briefly outline our second approach. We begin by bounding $E_f(X,S)$ in terms the quantity
    \[ \sum_{x,x' \in X} f(x) f(x') \sum_{s_1, \dots, s_{2t} \in S} \1_{x-x' = s_1 + \dots + s_t - s_{t+1} - \dots - s_{2t}} \]
	via an application of Lemma \ref{hldr}, where $t = 2m$. We call the latter quantity $E_t$. We divide such solutions into two types; the first of these having few \emph{cancelling pairs}, that is, pairs $(i,j) \in \{1,\dots, t\} \times \{1,\dots, t\}$ such that $s_i = s_j$. We upper bound the number of these solutions by crucially using the linear independence of the set $S$ along with the hypothesis that $(\sum_x f(x))^2$ is controlled nicely by $\sum_x f(x)^2$. The remaining type of solutions have many cancelling pairs. These we deal with an inductive, combinatorial argument; in particular, we show that the count of such solutions are controlled nicely by $E_{t-1}$. A second application of Lemma \ref{hldr} allows us to bound this in terms of $E_t$ and $\sum_x f(x)^2$. This combines with the preceding discussion to give us the desired bound on $E_f(X,S)$.

    Returning to Theorem \ref{thm:main}, our proof ideas are very analytic in nature. In particular, they crucially depend on Theorem \ref{asymenergy}. There is an alternative approach to proving Theorem \ref{thm:main} (found by a query to GPT-6 Astra; see \cref{comparisonremark} and \cref{appendix1}) which replaces the decomposition lemma and analytic estimates with a more geometric approach. The latter uses notions that are standard in convex geometry and relies crucially on the idea that, once a finite set $A$ is normalized appropriately, the intersection $A\cap(A+v)$ must be small if $\n{v}_2$ is large.
    (This is the content of \cref{lem:inter}.)

	\subsection*{Organization}
	In \S\ref{sec:prelim} we give some definitions and set some notation which we will use throughout the paper.
	We also state many of the ingredients we use to prove \cref{thm:main}.
	In \S\ref{sec:pf} we prove \cref{thm:main} assuming \cref{asymenergy}. 
	\S\S\ref{unweightedsection}--\ref{sec:hypercontractivity} contain three distinct proofs of statements related to \cref{asymenergy}. 
	First, in \S\ref{unweightedsection}, we prove \cref{hammingcor} using graph-theoretic techniques.
	In \S\ref{sec:asym}, we prove \cref{asymenergy} by analysing higher-order weighted additive energies.
	In  \S\ref{sec:hypercontractivity}, we provide an alternative proof of a weaker version of \cref{asymenergy} via hypercontractivity-type results. Finally, in Appendix~\ref{appendix1}, we give an alternative, more geometric proof of Theorem \ref{thm:main} which was found by GPT-6 Astra.

	\subsection*{Notation}
	In this paper we use Vinogradov notation; that is, we write $X \ll_z Y$, or equivalently $X = O_z(Y)$, to mean that $|X| \leq C_z Y$ for some constant $C_z>0$ which depends only on the parameter $z>0$. Given sets $A,B,C$, we write $A = B \sqcup C$ if $A = B \cup C$ and $B \cap C = \emptyset$. For any $n \in \mathbb{N}$, we denote $[n] = \{1,2,\dots, n\}$. For any finite, non-empty set $A \subseteq \R^d$, we write $\1_A$ for its characteristic function.

	\subsection*{Acknowledgements}
	We would like to thank Dan Altman for helpful conversations.  The first author is supported by a Leverhulme Early Career Fellowship ECF-2025-148. The second author is supported by the National Science Foundation Graduate Research Fellowship Program under Grant No.\ DGE-2146755. The first and third authors would like to thank the American Institute of Mathematics for providing excellent working conditions.

    \subsection*{AI Declaration}

    The main body of this paper, that is, the content in \S\S\ref{introductionsection}--\ref{sec:hypercontractivity}, was entirely generated and written by the authors without any AI assistance. AI tools, in particular GPT-6 Astra, were then used for literature search and proof-reading of this content. During this process, GPT-6 Astra suggested an alternative, more geometric approach to proving Theorem \ref{thm:main}. A significantly digested and rewritten version of this is recorded in Appendix \ref{appendix1}, with the authors having found various citations and sources for all the ideas used by GPT-6 Astra.

    \begin{remark} \label{comparisonremark}
    In this remark, we present a brief mathematical and philosophical comparison of the two proofs of Theorem \ref{thm:main}. The human-generated proof of Theorem \ref{thm:main}, that is, the proof recorded in \S\ref{sec:pf}, relies on heavy machinery from additive combinatorics (see Theorems \ref{thm:quantitative-BSG} and \ref{thm:weak-PFR}) along with various combinatorial and probabilistic manoeuvres and our human-generated proof of Theorem \ref{asymenergy}. The AI approach does not rely on the aforementioned heavy-machinery and does not consider the connection between Theorems \ref{thm:main} and \ref{asymenergy}. Instead, it builds a very direct geometric approach to studying Theorem \ref{thm:main}, relying on ideas from convex geometry to obtain a very explicit and optimal version of Theorem \ref{thm:main}, see Theorem \ref{thm:explicit}.

    Theorem \ref{thm:explicit} is a stronger version of Theorem \ref{thm:main} and the AI approach is more direct and shorter than the human-generated approach.  On the other hand, an advantage of the indirect human-generated approach is finding this connection between Theorem \ref{thm:main} and Theorem \ref{asymenergy}, which a priori appear to be somewhat unrelated. In fact, having to prove Theorem \ref{asymenergy} led us to find two new proofs of results of the shape \eqref{ksresult}, subsequently resulting in the improvements described in Theorem \ref{hammingcor}.  
    \end{remark}

	\section{Preliminaries}\label{sec:prelim}

	For two finitely-supported functions $f,g\colon\R^d\to\C$, we define their convolution $f*g\colon \R^d \to \C$ by setting
	\[(f*g)(x)=\sum_{y\in\supp f}f(y)g(x-y)\]
	for all $x \in \R^d$.
	For any $p > 0$ and any finitely-supported function $f\colon \R^d \to \mathbb{C}$, we define
	\[\n{f}_p = \Big(  \sum_{x \in \supp f} |f(x)|^p  \Big)^{1/p}  \ \ \text{and} \ \ \n{f}_\infty = \sup_x |f(x)| ,\]
	where $ \supp f = \{ x \in \mathbb{R}^d : f(x) \neq 0\}$.
	We will also write $E(X) = E(X,X)$ for any finite, non-empty set $X \subseteq \mathbb{R}^d$.

	Let $X,Y \subseteq \mathbb{Z}^d$ be finite, non-empty sets. Then we have the trivial bound
	\begin{equation}  \label{trivialbound}
		E(X,Y) \leq |X||Y|^2. 
	\end{equation}
	We have the following straightforward upper bounds on mixed additive energies. 
	
	\begin{lemma}\label{stp}
		Let $A,B,C,D$ be finite, non-empty sets in $\mathbb{R}^d$. Then
		\[ \sum_{a \in A}\sum_{b \in B}\sum_{c\in C}\sum_{d\in D} \1_{a+ b = c+d} \leq E(A,B)^{1/2} E(C,D)^{1/2}. \]
	\end{lemma}
	\begin{proof}  
		For any $x \in \mathbb{R}^d$, let
		\[ r(x) = \sum_{a \in A} \sum_{ b \in B} \1_{x = a+ b} \ \ \text{and} \ \ r'(x) = \sum_{c \in C} \sum_{ d \in D} \1_{x = c+ d} . \]
		Applying double-counting and Cauchy's inequality gives us 
		\[  \sum_{a \in A}\sum_{b \in B}\sum_{c\in C}\sum_{d\in D} \1_{a+ b = c+d} = \sum_{x \in A+B } r(x) r'(x) \leq \left(\sum_{x \in A+B} r(x)^2 \right)^{1/2} \left(\sum_{x \in A+B} r'(x)^2\right)^{1/2}. \]
		Using double-counting again, we see that
		\[ \sum_{x \in A+B} r(x)^2 = E(A,B) \ \ \text{and} \ \ \sum_{x \in A+B} r'(x)^2 \leq E(C,D) ,\]
		and so we are done.
	\end{proof}
	
	We will use the above lemma multiple times in our paper. For example, this immediately delivers the following result.
	\begin{lemma} \label{trivbd4}
		Let $X_1, X_2  \subseteq \mathbb{R}^d$ be finite, disjoint sets, let $X = X_1 \sqcup X_2$ and let $Y \subseteq \mathbb{R}^d$ be a finite, non-empty set. Then 
		\begin{equation} \label{trivialbound2} 
			E(X,Y) \leq 2E(X_1, Y) + 2E(X_2, Y). 
		\end{equation}
	\end{lemma}
	\begin{proof}
	
		Note that 
		\[ E(X,Y) = E(X_1, Y) + E(X_2, Y) + 2\sum_{x \in X_1} \sum_{x' \in X_2} \sum_{y_1, y_2 \in Y} \1_{x + y_1 = x' + y_2}. \]
		Lemma \ref{stp} and the AM--GM inequality imply that the last term on the right-hand side above can be bounded above by
		\[ 2 E(X_1,Y)^{1/2} E(X_2, Y)^{1/2} \leq E(X_1, Y) + E(X_2, Y). \]
		This concludes the proof of Lemma \ref{trivbd4}.
	\end{proof}
	
	\begin{remark} \label{alprem}
		It is worth noting that for any non-empty set $X \subseteq \{0,1\}^d$, Lemma \ref{stp} allows us to upper bound $T(X)$ in terms of $E(X,S)$ where $S = \{e_1, \dots, e_d\}$. Indeed, 
		\[ T(X) = \sum_{x_1, x_2 \in X} \sum_{1\leq i < j \leq d} (2 \cdot \1_{x_1 - x_2 = e_i - e_j} + 2\cdot\1_{x_1 - x_2 = e_i+ e_j}). \]
		Thus, we have
		\[ T(X) - (E(X,S) - |X||S|) = 2\sum_{x_1, x_2 \in X} \sum_{1 \leq i< j \leq d} \1_{x_1 - e_i = x_2 + e_j}=\sum_{x_1, x_2 \in X} \sum_{1 \leq i, j \leq d} \1_{x_1 - e_i = x_2 + e_j}. \]
		Lemma \ref{stp} implies that the right hand side is bounded above by $E(X,S)$, whence
		\[          E(X,S) - |X||S| \leq T(X) \leq 2E(X,S) - |X||S| . \]
	\end{remark}

	Let us now consider a slightly more general setting which arises in the proof of \cref{asymenergy}. Suppose $X, S$ are some finite, non-empty sets in $\mathbb{R}^d$ and $f\colon X \to [0, \infty)$ is a function and $k \geq 1$ is an integer. Define
	\begin{equation} \label{dn5}
		E_k = \sum_{x_1, x_2} \sum_{s_1, \dots, s_{2k} \in S} f(x_1) f(x_2) \1_{x_1 - x_2 = s_1 - s_2 + \dots + s_{2k-1} - s_{2k}}  \ \ \text{and} \ \ E_0 = \sum_{x \in X} f(x)^2.  
	\end{equation}
	We will often be interested in bounding $E_r$ in terms of $E_k$ for some $k \geq r \geq 0$.     When the set $X$ is contained in $\mathbb{Z}^d$,  this can be done in a straightforward manner by using orthogonality to pass to Fourier space and applying H\"{o}lder's inequality there. This can then be extended to sets $X \subset \R^d$ by using the theory of Freiman isomorphisms. This leads to the following result.

	\begin{lemma} \label{hldr}
		For any $r\geq 1$ and any $0 \leq k \leq r$, we have
		\[ E_k \leq E_0^{1 - k/r} E_r^{k/r}. \]
	\end{lemma}

Instead of using Fourier analysis, we will describe an alternative proof of Lemma \ref{hldr} which proceeds purely combinatorially.

    \begin{proof}[Proof of Lemma \ref{hldr}]
        We may assume $f$ does not vanish everywhere on $X$ since otherwise the result is trivial. This means that $E_0 = \sum_{x \in X} |f(x)|^2 >0$. Moreover, one can check that $E_k \geq \sum_{x \in X} |f(x)|^2 |S|^k > 0$ for every $k \geq 1$ by simply considering solutions satisfying $x_1 = x_2$ and $s_{2i-1} = s_{2i}$ for all $1 \leq i \leq k$. For any $j\in\Z_{\geq 0}$ and $y \in \R^d$, define
        \[ r_j(y) = \sum_{x\in X}\sum_{s_1, \dots, s_{j}\in S}f(x) \1_{y = x + s_1 + \dots + s_j} \ \ \text{and} \ \ r_j'(y) = \sum_{x\in X}\sum_{s_1, \dots, s_{j}\in S}f(x) \1_{y = x + s_1 + \dots + s_{j-1}-  s_j} \]
        For any positive integer $j$, double-counting gives us
        \[ E_j = \sum_{y\in \R^d}r_j(y)^2=\sum_{y \in \R^d} r_j'(y)^2  = \sum_{z \in \R^d} r_{j+1}'(z) r_{j-1}(z). \]
        Applying Cauchy's inequality and double-counting, we get that
        \begin{equation}\label{eq:E-convex}
            E_j \leq E_{j+1}^{1/2} E_{j-1}^{1/2}.
        \end{equation}

        Now, for any $0 \leq j \leq r$, let $a_j = \log E_j$ and for any $1 \leq j \leq r$, let $\Delta_j = a_j - a_{j+1}$. The bound \eqref{eq:E-convex} is equivalent to writing $\Delta_{j} \leq \Delta_{j+1}$ for any $1 \leq j \leq r$. This means that for any $1 \leq k \leq r$, one has
        \[ k^{-1}\sum_{i=1}^k \Delta_i \leq r^{-1} \sum_{i=1}^r \Delta_i.  \]
        Simplifying, we get that $(a_k - a_0)/k \leq (a_r - a_0)/r$, that is, 
        \[ a_k \leq a_r (k/r) + a_0 (1 - k/r).\]
        This is precisely the desired bound.
    \end{proof}

	\begin{remark}
    This proof is reminiscent of the ideas described in \cite{Mu2026}, and in fact, Lemma \ref{hldr} follows from   \cite[Theorem 1.1]{Mu2026}. Variations of this type of idea, when $k =1$ and $r= 2^l$ for some $l \in \mathbb{N}$, have also been used in the proof of the asymmetric Balog--Szemer\'{e}di--Gowers theorem, see \cite[Lemma 42]{Shk2024}.
	\end{remark}

	In our proof of \cref{thm:main}, we will find it useful to pass from an upper bound on the additive energy between two sets $X$ and $Y$ to an upper bound on $E(X,Y')$ where $Y'$ is a small set sampled from $Y$.
	If $Y'$ is obtained by taking $s$ elements of $Y$ uniformly at random, it is not too hard to upper-bound the expected value of $E(X,Y')$ in terms of $E(X,Y)$.
	We require a slight generalization of such a statement where the elements of $Y'$ are not too correlated.
	
	\begin{lemma}[Subsampling]\label{lem:subsample}
		Let $\mathcal D$ be a probability distribution on tuples $(y_1,\ldots,y_s)$ of distinct elements of $Y$.
		Suppose that, for each $1\leq i<j\leq s$, the marginal distribution of the pair $(y_i,y_j)$ has maximum point probability at most $C/\abs{Y}^2$ for some parameter $C>0$.
		Then
		\[\EE_{(y_1,\ldots,y_s)\sim\mathcal D} E(X,\{y_1,\ldots,y_s\})\leq s\abs{X}+\frac{Cs^2}{\abs{Y}^2}E(X,Y).\]
	\end{lemma}
	\begin{proof}
		We have
		\[E(X,\{y_1,\ldots,y_s\})=s\abs{X}+2\sum_{1\leq i<j\leq s}(\1_X*\1_{-X})(y_i-y_j).\]
		Our assumption on the point probabilities of pairs $(y_i,y_j)$ implies that, for $i\neq j$,
		\begin{align*}
			\EE_{(y_1,\ldots,y_s)\sim\mathcal D}(\1_X*\1_{-X})(y_i-y_j)
			&=\sum_{z_1,z_2\in Y}\Pr[(y_i,y_j)=(z_1,z_2)](\1_X*\1_{-X})(z_1-z_2)\\
			&\leq \frac C{\abs{Y}^2}\sum_{z_1,z_2\in Y}(\1_X*\1_{-X})(z_1-z_2)=\frac C{\abs{Y}^2}E(X,Y).
		\end{align*}
		The claimed bound follows.
	\end{proof}

	
	We will need the following quantitative version of the Balog--Szemer\'edi--Gowers theorem due to Reiher--Schoen \cite[Theorem 1.2]{RS2024}.
	For our application, the quality of the upper bound on $\abs{A'-A'}$ is not so important, but it is crucial to have a bound of the form $\abs{A'}\gg\delta^c\abs{A}$ for some $c<1$.
	
	\begin{theorem} \label{thm:quantitative-BSG}
		Let $A$ be a finite nonempty subset of an abelian group and let $\delta>0$. 
		Suppose that $E(A,A)\geq \delta\abs{A}^3$. 
		Then there exists some subset $A'\subseteq A$ for which 
		\[ \abs{A'}\gg\delta^{1/2}\abs{A} \ \ \text{and} \ \ \abs{A'-A'}\ll\delta^{-4}\abs{A'}. \]
	\end{theorem}

	We will also need a structure theorem for sets with small doubling in $\R^d$. In particular, the following can be deduced by combining the resolution of the weak polynomial Freiman--Ruzsa conjecture due to Gowers--Green--Manners--Tao \cite[Theorem~1.3]{GGMT2025} along with a standard application of the Pl\"{u}nnecke--Ruzsa inequality and Ruzsa's covering lemma.

	\begin{theorem} \label{thm:weak-PFR} 
		There exists an absolute constant $C > 0$ such that for every $d \in \mathbb{N}$, any finite, non-empty $A\subset\R^d$ satisfying $\abs{A-A}\leq K\abs{A}$, for some $K \geq 1$, is contained within $K^C$ translates of some vector space of dimension at most $C\log K$.
	\end{theorem}
	\begin{proof}
		By the Pl\"unnecke--Ruzsa inequality, we have $\abs{A+A}\leq K^2\abs{A}$.
		The abelian group $\langle A\rangle$ generated by $A$ is finitely generated and torsion-free, and so it is isomorphic to $\Z^n$ for some positive integer $n$. 
		\cite[Theorem~1.3]{GGMT2025} implies that there is some $A'\subset A$ with $\abs{A'}\geq K^{-C_1}|A|$ and with affine span in $\Z^n$ of dimension at most $2C_2\log K$, where $C_1$ and $C_2$ are absolute constants (the same as those in \cite[Theorem~1.3]{GGMT2025}).
		
		We claim that the dimension of the affine span of $A'$ in $\R^d$ is at most that in $\Z^n$.
		Indeed, the statement that $A'$ has affine span in $\Z^n$ of dimension $D\leq 2C_2\log K$ implies that, if $a_1,\ldots,a_{D+2}\in A'$, then there exist $c_1,\ldots,c_{D+2}\in\Q$, not all zero but summing to zero, for which 
		\begin{equation}\label{eq:lindep}
			c_1a_1+\cdots+c_{D+2}a_{D+2}=0
		\end{equation}
		Scaling up, we can take $c_1,\ldots,c_{D+2}\in\Z$ and interpret \eqref{eq:lindep} in $\R^d$. From here, we obtain that the affine span of $A'$ in $\R^d$ has dimension at most $D$.
		
		Finally, let $V$ be the subspace of $\R^d$ spanned by $A'-A'$.
		The desired result now follows from Ruzsa's covering lemma: since $\abs{A+A'}\leq \abs{A+A}\leq K^2\abs{A}\leq K^{C_1+2}\abs{A'}$, the set $A$ can be covered by at most $K^{C_1+2}$ translates of $A'-A'\subset V$.
	\end{proof}
	
	\section{Proof of Theorem \ref{thm:main}}\label{sec:pf}
	
	In this section, we prove \cref{thm:main} assuming \cref{asymenergy}.
	The first step is to find a decomposition of our set $A$ into two pieces, one of which is structured and the other of which is either small or not too structured. For this purpose, we will find the following notation useful: given a non-empty set $S \subseteq \mathbb{R}^d$ and some subspace $V$ of $\mathbb{R}^d$, we denote by $S/V$ the image $\pi(S)$ of $S$ under the canonical projection map $\pi\colon\mathbb{R}^d \to \mathbb{R}^d/V$.

	\begin{lemma}\label{lem:decomp}
		Suppose $d$ is sufficiently large. Set $\eps:=d^{-1/10}$ and $\delta:=d^{-3/2}$.
		Then any finite set $A\subset\R^d$ with $\dim(A) = d$ can be decomposed into two sets $A_{\mathrm{str}}\sqcup A_{\mathrm{unstr}}$ such that the following two properties hold:
		\begin{enumerate}
			\item The set $A_{\mathrm{str}}$ is contained within the union of $d^{O(1)}$-many translates of some subspace $V\subset\R^d$ of dimension at most $d/2$.
			
			\item The set $A_{\mathrm{unstr}}$ either satisfies $\abs{A_{\mathrm{unstr}}}\leq\frac1d\abs{A}$, or else satisfies both (i) that $E(A_{\mathrm{unstr}})<\delta\abs{A_{\mathrm{unstr}}}^3$ and (ii) that, for any subspace $V\subset V'\subset\R^d$ with $\dim (V'/V) \leq d/100$ and any $x\in\R^d$, we have
			\begin{equation} \label{assum678}
				\abs{A_{\mathrm{unstr}}\cap (x+V')}<(1-\eps)\abs{A_{\mathrm{unstr}}}.
			\end{equation}
		\end{enumerate}
	\end{lemma}
	\begin{proof}
		We will construct our decomposition iteratively.    Initially, set $A_1=\emptyset$ and $B_1=A$ and $V_1 = \{ 0\}$ and $\rho_1 = |B_1|/|A| =1$. At the beginning of the $i^{\rm th}$ step of the procedure, for any integer $i \geq 1$, we will have disjoint sets $A_i$ and $B_i$ such that $A = A_i \sqcup B_i$, the parameter $\rho_i = |B_i|/|A|$ and the vector space $V_i \subset \mathbb{R}^d$. With this in hand, the $i^{\rm th}$ step involves proceeding through the following conditions in order. 
		
		\begin{enumerate}[1.]
			\item If $\rho_i\leq1/d$, terminate the procedure. \label{step1}
			
			\item If there exists a subset $B' 
			\subseteq B_i$ such that $|B'| \geq (1-\eps)|B_i|$ and $B'/V_i \subseteq (x+V')/V_i$ for some $x \in \R^d$ and some subspace $V_i\subset V'\subset\R^d$ satisfying $\dim V'-\dim V_i\leq d/100$, update
			\begin{align*}
				V_{i+1}:=V' \ \ \text{and} \ \  A_{i+1}:=A_i\sqcup (B_i\cap (x+V')) \  \ \text{and}  \ \            B_{i+1}:=B_i\setminus (B_i\cap (x+V')).
			\end{align*}
			In this way, we have $\rho_{i+1}\leq \eps\rho_i$. Proceed to the $(i+1)^{\rm th}$ step. \label{step2}
			
			\item If instead $E(B_i,B_i)\geq \delta\abs{B_i}^3$, apply \cref{thm:quantitative-BSG} and then \cref{thm:weak-PFR} to find some subset $T\subset B_i$ for which $\abs{T}\geq c\delta^{1/2}\abs{B_i}$ and which may be covered by $\delta^{-C'}$ translates of some vector space $W$ of dimension at most $C'\log(1/\delta)$, where $c,C'>0$ are constants. Update
			\begin{align*}
				V_{i+1}=V_i+W \ \ \text{and} \ \  A_{i+1}=A_i\sqcup T \ \ \text{and} \ \ 
				B_{i+1}=B_i\setminus T.
			\end{align*}
			In this way, we have $\rho_{i+1}\leq(1-c\delta^{1/2})\rho_i$. Proceed to the $(i+1)^{\rm th}$ step. \label{step3}
			
			\item If neither of the previous steps have been successful, terminate the procedure.  \label{step4}
		\end{enumerate}

		It is simple to check that this procedure can only run for finitely many steps, and so, suppose the above procedure stops after $t$ steps, for some $t \in \mathbb{N}$. The termination condition in event \ref{step1} and our choices of $\eps$ and $\delta$ imply that event \ref{step2} can occur at most $10$ times and event \ref{step3} can occur at most $c^{-1}d^{3/4}\log d$ times. We conclude that, at the end of the procedure,
		\[\dim V_t\leq \frac d{10}+O(d^{3/4}(\log d)^2)\leq \frac d2.\]
		Moreover, at the $i^{\rm th}$ step of the procedure, we observe that $A_i$ is covered by at most $i\cdot d^{C''}$ translates of $V_i$ for some absolute constant $C''$.
		We thus conclude that $A_t$ is covered by at most $d^{O(1)}$ translates of $V_t$.
		We now set $A_{\mathrm{str}}:=A_t$ and $A_{\mathrm{unstr}}:=B_t$ and $V:=V_t$ to obtain the desired conclusion.
	\end{proof}
	
	Next, using \cref{asymenergy}, we bound the additive energy between $A_{\mathrm{str}}$ and any set which remains affine-independent when quotiented by $V$. 
	We will in fact need the following slight generalization of \cref{asymenergy}.
	
	\begin{lemma}\label{lem:asym-affind}
		Let $d,m\in\N$ with $m\leq d$.
		Let $S$ be an affine-independent subset of $\R^d$, and let $f\colon\R^d\to[0,\infty)$ be supported on a finite, non-empty set $X\subseteq\R^d$ and satisfy $(\n{f}_1/\n{f}_2)^2\leq\binom dm$.
		Then we have
		\[E_f(X,S)\ll md\n{f}_2^2.\]
	\end{lemma}
	\begin{proof}
		Since translating $S$ does not change $E_f(X,S)$, we may translate $S$ suitably so as to ensure that $0\in S$. By changing the coordinate frame, we may assume that the other $\abs{S}-1$ elements of $S$ are $e_1,\ldots,e_j$ for $j=\abs{S}-1$.
		We conclude from \cref{asymenergy} that
		\begin{align*}
			E_f(X,S)
			&\leq E_f(X,\{0,e_1,\ldots,e_d\})\\
			& \leq \sum_{x\in X}f(x)^2+2\sum_{i=1}^d\sum_{x\in X}f(x)f(x+e_i)+E_f(X,\{e_1,\ldots,e_d\})\\
			&\ll (d+1+md)\n{f}_2^2.\qedhere
		\end{align*}
	\end{proof}
	
	We will use this lemma to show that any set $S \subset \R^d$ with $|S| =|S/V| \gg d$ and with $S/V$ being affine-independent  has little additive interaction with $A_{\rm str}$. 
	
	\begin{lemma}\label{lem:structure}
		Let $V$ and $A_{\mathrm{str}}$ be as in the conclusion of \cref{lem:decomp}.
		Let $S$ be any set of size at least $d/100$ such that $S/V$ is affine-independent and $|S| = |S/V|$.
		Then $E(A_{\mathrm{str}},S)\ll\abs{A}\abs{S}$.
	\end{lemma}
	\begin{proof}
		Let $\pi\colon \mathbb{R}^d \to \mathbb{R}^d/V$ be the canonical projection map, let $X = \pi(A_{\mathrm{str}})$ and let $S' = \pi(S)$.
		Condition (1) of \cref{lem:decomp} implies that $\abs{X}\leq d^{O(1)}$. For any $x \in X$, let $A_x = A_{\mathrm{str}} \cap (x + V)$. Note that
		\begin{align*}
			E(A_{\mathrm{str}},S) & = \sum_{x_1, x_2 \in X}  \sum_{s_1, s_2 \in S}   |A_{x_1} \cap (A_{x_2} + s_1 - s_2)| \\
			&  \leq  \sum_{x_1, x_2 \in X}  \sum_{s_1, s_2 \in S} \1_{x_1 - x_2 = \pi(s_1) - \pi(s_2)} |A_{x_1}|^{1/2} | A_{x_2}|^{1/2} .
		\end{align*}
		In the second step above, we have used the fact that whenever the set $A_{x_1} \cap (A_{x_2} + s_1 - s_2)$ is non-empty, we must have $x_1 = x_2 + \pi(s_1) - \pi(s_2)$.  Define $f\colon X \to [0, \infty)$ by setting $f(x)=|A_x|^{1/2}$ for every $x\in X$.
        We have the bound
        \[\frac{\n{f}^2_1}{\n{f}_2^2}\leq\abs{X}\leq d^{O(1)}\leq(\dim(\R^d/V))^{O(1)}.\]
        So, we may apply \cref{lem:asym-affind} to $X,S'\subset\R^d/V$ with $m=O(1)$ to obtain
		\[E(A_{\mathrm{str}},S)\leq E_f(X,S')\ll \abs{A_{\mathrm{str}}}\abs{S'}=\abs{A_{\mathrm{str}}}\abs{S}\leq\abs{A}\abs{S}.\qedhere\]
	\end{proof}
	
	Now, we provide a procedure for selecting such a set $S$ for which $E(A_{\mathrm{unstr}},S)$ is well-controlled.
	
	\begin{lemma}\label{lem:unstructure}
		Let $V$ and $A_{\mathrm{unstr}}$ be as in \cref{lem:decomp}.
		Then there exists a subset $S$ of $A$ of size at least $d/100$ such that $S/V$ is affine-independent and $|S| = |S/V|$  and \[ E(A_{\mathrm{unstr}},S)\ll\abs{A}\abs{S}. \]
	\end{lemma}
	\begin{proof}
		We first remark that $A/V$ is full-dimensional in $\R^d/V$, and so there exists a subset $S \subseteq A$ such that $|S| \geq d/100$ and $S/V$ is affine-independent and $|S| = |S/V|$. If $\abs{A_{\mathrm{unstr}}}\leq\frac1d\abs{A}$, then the trivial bound \eqref{trivialbound} implies
		\[E(A_{\mathrm{unstr}},S)\leq\abs{A_{\mathrm{unstr}}}\abs{S}^2\ll\abs{A}\abs{S}.\]
		So, we may henceforth assume that we are in case (2) of \cref{lem:decomp}.
		
		In this case, we use a randomized procedure to construct such a set $S$, with an aim of applying \cref{lem:subsample}.    For $1\leq j\leq d/100$, choose an element $y_j\in A_{\mathrm{unstr}}$ uniformly at random subject to the condition that the set $\{y_1,\ldots,y_j\}/V$ is affine-independent.
		There are at least $\eps\abs{A_{\mathrm{unstr}}}$ choices for each $y_j$: otherwise, the affine-span $V'$ of $\{y_1,\ldots,y_{j-1}\}+V$ would satisfy $|V' \cap A_{\rm unstr}| \geq (1 - \varepsilon)|A_{\rm unstr}|$, contradicting \eqref{assum678}.
        In particular, given the values of $y_1,\ldots,y_{j-1}$, the maximum point probability of $y_j$ is at most $\eps^{-1}/\abs{A_{\mathrm{unstr}}}$.
		We claim moreover that, for each $1\leq i<j\leq d/100$, the maximum point probability of the pair $(y_i,y_j)$ is at most $\eps^{-2}/\abs{A_{\mathrm{unstr}}}^2$.
        Indeed, for each $a_1,a_2\in A$,
        \begin{align*}
            \Pr[y_i=a_1,y_j=a_2]
            &\leq\max_{z_1,\ldots,z_{i-1}}\Pr\big[y_i=a_1,y_j=a_2\mid y_1=z_1,\ldots,y_{i-1}=z_{i-1}\big]\\
            &=\max_{z_1,\ldots,z_{i-1}}\Pr\big[y_i=a_1\mid y_{[1,i)}=z_{[1,i)}\big]\cdot\Pr\big[y_j=a_2\mid y_{[1,i)}=z_{[1,i)},y_i=a_1\big]\\
            &\leq\frac{\eps^{-1}}{\abs{A_{\mathrm{unstr}}}}\cdot\max_{z_1,\ldots,z_{i-1}}\Pr\big[y_j=a_2\mid y_{[1,i)}=z_{[1,i)},y_i=a_1\big]\\
            &\leq\frac{\eps^{-1}}{\abs{A_{\mathrm{unstr}}}}\cdot\max_{z_1,\ldots,z_{j-1}}\Pr\big[y_j=a_2\mid y_{[1,i)}=z_{[1,i)},y_i=a_1,y_{(i,j)}=z_{(i,j)}\big]\\
            &\leq\frac{\eps^{-2}}{\abs{A_{\mathrm{unstr}}}^2},
        \end{align*}
        where we have used notation like $y_{[1,i)}$ as shorthand for the tuple $(y_1,\ldots,y_{i-1})$.
		Therefore, letting $S=\{y_1,\ldots,y_{\lfloor d/100\rfloor}\}$, we may apply  \cref{lem:subsample} to deduce that
		\[\EE E(A_{\mathrm{unstr}},S)\leq \abs{A_{\mathrm{unstr}}}\abs{S}+\frac{\eps^{-2}(d/100)^2}{\abs{A_{\mathrm{unstr}}}^2}E(A_{\mathrm{unstr}})\ll\abs{A}\abs{S}+\eps^{-2}d^2\delta\abs{A_{\mathrm{unstr}}}\ll\abs{A}\abs{S}.\]
		This means that there must exist some set $S \subseteq A_{\rm unstr}$ satisfying the conclusion of \cref{lem:unstructure}, and so, we are done.
	\end{proof}
	
	\cref{thm:main} follows from combining Lemmas \ref{lem:decomp}, \ref{lem:structure} and \ref{lem:unstructure} with \eqref{trivialbound2}.

	\section{The edge-isoperimetric problem over the Boolean cube} \label{unweightedsection}
	
	In this section, we present our proof of Theorem  \ref{hammingcor}. 
	So that we can use additive notation, we identify the vertices of the hypercube with $(\Z/2\Z)^d$. To avoid ambiguity with $\Z^d$ or $\R^d$, we write $e_i'\in(\Z/2\Z)^d$ for the vector with a $1$ in coordinate $i$ and zeros in other coordinates. In particular, $\{e_1',\dots, e_d'\}$ forms a basis for $(\mathbb{Z}/2\mathbb{Z})^d$. Moreover, two vertices $x_1,x_2\in(\Z/2\Z)^d$ have Hamming distance $2$ if and only if their difference $x_1-x_2$ can be written as $e_i'+e_j'$ for some distinct $1\leq i,j\leq d$. This observation immediately delivers the  bound 
	\begin{equation} \label{atrivbd}
		T(X) \leq |X| { d \choose 2} . 
	\end{equation}
	We are now ready to present our proof of \cref{hammingcor}.
	
	\begin{proof}[Proof of \cref{hammingcor}]
		Let $G$ be the graph with vertex set $X$ where vertices $x_1$ and $x_2$ are joined by an edge if $x_1-x_2=e_i'+e_j'$ for some distinct $i,j\in\{1,\ldots,d\}$. 
		Two edges with increments $e_i'+e_j'$ and $e_k'+e_\ell'$ are said to \emph{conflict} if $\{i,j\}\cap\{k,\ell\}$ is nonempty. 
		We call a subset of the edge set of $G$ \emph{non-conflicting} if no two edges in the subset conflict.
		The number of edges of $G$ is $T(X)/2$, so the average degree $L$ of $G$ satisfies $L=T(X)/\abs{X}$.
		We may assume for the sake of contradiction that $L\geq 4md$.
		
		Let $X'\subset X$ be a subset such that the induced subgraph $G':=G[X']$ has minimum degree exceeding $L/2$. Such a subset can be found by repeatedly removing vertices of degree at most $L/2$ until none remain; each removal does not decrease the average degree, so this process terminates in a nonempty graph.
		Let $P$ be the number of paths of length $m$ consisting of non-conflicting edges in $G'$. Our desired bound will follow from obtaining suitable upper and lower bounds for $P$.
		
		For the upper bound, note that between any two distinct vertices $x$ and $y$ in $G'$, there are at most $(2m)!/2^m$ non-conflicting paths.
		Indeed, if there is any such path, then we can write
		\[x-y=e_{i_1}'+e_{j_1}'+\cdots+e_{i_m}'+e_{j_m}'\]
		for some distinct $i_1,j_1,\ldots,i_m,j_m\in\{1,\ldots,d\}$. Since the set $\{e_1',\ldots,e_d'\}$ is linearly independent, such a representation is unique up to permutation of $\{i_1,j_1,\ldots,i_m,j_m\}$. The factor of $2^{-m}$ comes from the fact that swapping $i_k$ and $j_k$ does not change the path. We thus conclude that
		\begin{equation}\label{eq:P-upper}
			P\leq\frac{(2m)!}{2^m}\abs{X'}^2.
		\end{equation}
		
		On the other hand, given a non-conflicting path consisting of edges 
		\[ \{(x_0,x_1),(x_1, x_2),\ldots,(x_{t-1},x_t)\}, \]
		there are at most $2dt$ edges out of $x_t$ which conflict with some edge previously in the path. 
		Indeed, at most $2d$ putative new edges $e_k'+e_\ell'$ conflict with each previous $e_{i_r}'+e_{j_r}'$, as at most $2d$ subsets $\{k,\ell\}\subset\{1,\ldots,d\}$ of size $2$ intersect $\{i_r,j_r\}$. Thus, given any $x\in X'$, there are at least
		\[\prod_{r=0}^{m-1}\left(\frac L2-2dr\right)\geq\prod_{r=0}^{m-1}(2d(m-r))=(2d)^mm!\]
		non-conflicting paths starting from the vertex $x$. We conclude
		\begin{equation}\label{eq:P-lower}
			P\geq(2d)^mm!\abs{X'}.
		\end{equation}
		Combining \eqref{eq:P-upper} and \eqref{eq:P-lower} gives
		\[\abs{X'}\geq \frac{2^m(2d)^mm!}{(2m)!}=\frac{4^m}{\binom{2m}m}\cdot\frac{d^m}{m!}>\binom dm.\]
		This contradicts the fact that $\abs{X'}\leq\abs{X}\leq\binom dm$, and so we are done.
	\end{proof}
	
	\begin{remark}\label{rmk:strength}
		For $m\gg d$, our proof of \cref{hammingcor} can be pushed slightly further. 
		Namely, the upper bound $2dt$ on the number of potential conflicting edges can be strengthened to $\binom d2-\binom{d-2t}2=t(2d-2t-1)$, which is a substantial improvement when $t\gg d$. Following this, one can prove that $\abs{X}\leq\binom{d-m-1}m$ implies $T(X)\leq4m(d-m)\abs{X}$, which improves the leading constant in \cref{hammingcor} for some ranges of $m$.
	\end{remark}
	
	\begin{remark}\label{rmk:gen}
		The proof of \cref{hammingcor} we just presented can be modified to prove a bound
		\[E(X,S)\leq (4m+1)d\abs{X}\]
		for $X\subset\R^d$ of size at most $\binom dm$, where $S=\{e_1,\ldots,e_d\}$.
		Construct the graph $G$ with vertex set $X$ and where $x_1,x_2\in X$ are joined by an edge if $x_1-x_2=e_i-e_j$ for some distinct $i,j\in\{1,\ldots,d\}$. The average degree  of $G$ is $L:=E(X,S)/\abs{X}-d$. Call two edges $e_i-e_j$ and $e_{i'}-e_{j'}$ of $G$ \emph{conflicting} if $i=j'$ or $i'=j$. 
		Assuming $L\geq 4md$ and following the same strategy, the bounds \eqref{eq:P-upper} and \eqref{eq:P-lower} read
		\[(2d)^mm!\abs{X'}\leq P\leq m!^2\abs{X'}^2,\]
		and the bound $\abs{X'}>\binom dm$ can be extracted similarly.
	\end{remark}

	\section{Weighted Asymmetric Additive Energy Estimates} \label{sec:asym}
	
	In this section, we present our proof of Theorem \ref{asymenergy}.
	This requires an analysis of higher-order energies. So, for each integer $s \geq 1$ and each $y\in\R^d$, we define the quantities
	\[
	\mathcal E_s(y):=\{(i_1,\ldots,i_s,j_1,\ldots,j_s)\in[d]^{2s}:y=e_{i_1}-e_{j_1}+\cdots+e_{i_s}-e_{j_s}\} \]
	and
	\[
	\mathcal R_s(y) :=\{(i_1,\ldots,i_s,j_1,\ldots,j_s)\in\mathcal E_s(y):i_k\neq j_\ell\text{ for all }1\leq k,\ell\leq s\}.\]
	Moreover, let $\mathcal R_0(0) = \mathcal E_0(0)  = \{1\}$ and $\mathcal R_0(y) = \mathcal E_0(y) = \emptyset$ for all $y \neq 0$. Now, for any integer $s \geq 0$,  write
	\[ R_s = \sum_{x_1, x_2\in X}f(x_1)f(x_2)\abs{\mathcal R_s(x_1-x_2)} \]
	and
	\begin{align} \label{db2}
		E_s & = \sum_{x_1, x_2\in X}f(x_1)f(x_2)\abs{\mathcal E_s(x_1-x_2)} \nonumber \\
		& =\sum_{x_1, x_2 \in X} f(x_1) f(x_2)\sum_{\substack{1\leq i_1,\ldots,i_s\leq d\\1\leq j_1,\ldots,j_s\leq d}} \1_{x_1 - x_2 = e_{i_1}-e_{j_1}+\cdots+e_{i_s}-e_{j_s}}. 
	\end{align}
	Note that the last of these definitions is exactly the one recorded in \eqref{dn5}. In particular,
	\begin{equation}  \label{obs5}
		E_0=\n{f}_2^2 \ \ \text{and} \ \ E_1=E_f(X,S).
	\end{equation}
    We must show that $E_1\leq 4d(m+1/m)E_0$ subject to the condition that $\n{f}_1^2/\n{f}_2^2\leq\binom dm$.
    We may assume $m\leq d/2$.
	
	The main estimate behind the proof of Theorem \ref{asymenergy} is as follows.
    
	\begin{lemma} \label{claim4}
		Let $t = 2m$. Then
		\[ E_t\leq 2dtE_{t-1}+(2d)^tt!E_0. \]
	\end{lemma}
	
	This combines naturally with Lemma \ref{hldr} to imply the conclusion of Theorem \ref{asymenergy}.
	
	\begin{proof}[Proof of Theorem \ref{asymenergy}]
		First, suppose that
		\begin{equation} \label{case41}
			(2d)^t t! E_0 \leq  (2dt)E_{t-1}/m^2.
		\end{equation} 
		In this case, Lemma \ref{claim4} implies that
		\[ E_t \leq (2dt)E_{t-1}( 1+ 1/m^2),\]
		whence Lemma \ref{hldr} implies that
		\[ E_t^{1/t}  \leq (2dt)(1 + 1/m^2) E_0^{1/t}. \]
		A second application of Lemma \ref{hldr} then gives us
		\[ E_1 \leq E_0^{1-1/t} E_t^{1/t} \leq (2dt)(1 + 1/m^2) E_0. \]
		Combining this with \eqref{obs5} yields the required bound. Now, if \eqref{case41} does not hold, then we get that
		\[ E_{t-1} < m^2 (2d)^{t-1} (t-1)! E_0, \]
		whence Lemma \ref{hldr} gives us
		\[ E_1 \leq E_{t-1}^{1/(t-1)} E_0^{1 - 1/(t-1)} < m^{2/(t-1)} (2d) ((t-1)!)^{1/(t-1)} E_0.  \]
		Since $t= 2m$, we can use the AM-GM inequality to get $((t-1)!)^{1/(t-1)} \leq m$.  Moreover, one can see that $m^{1/(m-1/2)} \leq 2$ for all integers $m \geq 1$, which dispenses the required estimate in this case.
	\end{proof}
	
	Thus we will now prepare to prove Lemma \ref{claim4}. We will start by analysing $\mathcal E_s(y)$ for fixed choices of $y$ and integer $s \geq 2$. In particular, for any tuple $\tau = (i_1, \dots, i_s, j_1, \dots, j_s) \in \mathcal E_s(y)$, we call the pair $(i_k, j_\ell)$ for some $k, \ell \in [s]$ a \emph{cancelling pair} if $e_{i_k} - e_{j_\ell} = 0$. Furthermore, write $h(\tau)$ for the maximum number of disjoint cancelling pairs in $\tau$. For any $0 \leq u \leq s$, let 
	\[ \mathcal E_s^u(y) = \{ \tau \in \mathcal E_s(y) : h(\tau) = u\}. \]
	Note that for any $\tau \in \mathcal E_s^u(y)$ and any $0 \leq u < s$, once we remove the $u$ disjoint cancelling pairs, we obtain an element of $\mathcal R_{s-u}(y)$. This gives us the bound
	\begin{equation} \label{bound7}
		|\mathcal{E}_s^u(y)| \leq  {s \choose u}^2 d^u u!|\mathcal R_{s-u}(y)|,
	\end{equation}
	since we can choose the indices corresponding to the $u$ cancelling pairs in ${s \choose u}^2$ many ways and fix the values of these indices in $d^u u!$ many ways.

	With this in hand, we record a useful intermediate bound. 
	
	\begin{lemma}  \label{claim57}
		For any $1 \leq u \leq s$, we have
		\[ |\mathcal{E}_s^u(y)| \leq d s^2 |\mathcal{E}_{s-1}^{u-1}(y)|/u. \]
	\end{lemma}
	\begin{proof}
        Construct a bipartite \emph{multigraph} with the vertex set $\mathcal{E}_s^u(y) \sqcup \mathcal{E}_{s-1}^{u-1}(y)$ by connecting any $\tau \in \mathcal{E}_{s}^u(y)$ to $\tau' \in \mathcal{E}_{s-1}^{u-1}(y)$ if one obtains $\tau'$ from $\tau$ by removing one cancelling pair. This graph is a multigraph since we are allowing for multiple edges between two vertices; indeed, some $\tau$ may have multiple choices of cancelling pairs which, when cancelled, give us the same $\tau' \in \mathcal{E}_{s-1}^{u-1}(y)$. Let $E'$ be the number of edges in this graph. Note that
		\[  E' \geq u |\mathcal{E}_s^u(y)| \]
		since every $\tau \in \mathcal{E}_s^u(y)$ has at least $u$ cancelling pairs, each of which, when removed, give us some element in $\mathcal{E}_{s-1}^{u-1}(y)$. On the other hand, we have
		\[ E' \leq d s^2 |\mathcal{E}_{s-1}^{u-1}(y)|\]
		since for any fixed $\tau \in \mathcal{E}_{s-1}^{u-1}(y)$, we have at most $s^2$ choices for the positions where we may insert the two extra indices associated to a cancelling pair and at most $d$ choices for the value that this cancelling pair is to take. Combining the above two bounds delivers the desired bound, thus finishing our proof of Lemma \ref{claim57}.
	\end{proof}

	We are now ready to prove Lemma \ref{claim4}.
	
	\begin{proof}[Proof of Lemma \ref{claim4}]
		Using \eqref{db2} and the preceding discussion, we see that
		\[ E_t = \sum_{x_1, x_2 \in X} f(x_1) f(x_2) | \mathcal{E}_t(x_1 - x_2)|  =  \sum_{x_1, x_2 \in X} f(x_1) f(x_2)   \sum_{u=0}^t |\mathcal{E}_t^u(x_1 - x_2)|. \]
		Using \eqref{bound7} and Lemma \ref{claim57}  to estimate the cases when $0 \leq u < t/2$ and $t/2 \leq u \leq t$ respectively, we get that
		\begin{equation} \label{db3}
			E_t \leq \sum_{0 \leq u < t/2} {t \choose u}^2 d^u u! R_{t-u} + \sum_{x_1, x_2 \in X} f(x_1) f(x_2) \sum_{t/2 \leq u \leq t} \frac{dt^2}u | \mathcal{E}_{t-1}^{u-1}(x_1 - x_2)|.
		\end{equation}
		The second sum on the right hand side of \eqref{db3} can be bounded above by
		\begin{equation} \label{db4}
			2dt \sum_{x_1, x_2 \in X} f(x_1) f(x_2) \sum_{t/2 \leq u \leq t} |\mathcal{E}_{t-1}^{u-1}(x_1 - x_2)| \leq 2dt E_{t-1}. 
		\end{equation}
		We now consider the first sum, and so let $0 \leq u < t/2$. For any $x_1,x_2\in X$, a representation of $x_1-x_2$ as $e_{i_1}-e_{j_1}+\cdots+e_{i_{t-u}}-e_{j_{t-u}}$ with the $i_k$ distinct from the $j_\ell$, if such a representation exists, is unique up to permutation of the $i_k$ and the $j_\ell$. This follows from the linear independence of the set $\{e_1, \dots, e_d\}$.   We conclude that
		\[R_{t-u}\leq((t-u)!)^2\sum_{x_1,x_2\in X}f(x_1)f(x_2) = ((t-u)!)^2\n{f}_1^2.\]
		Our hypothesis implies that $\n{f}_1^2 \leq \n{f}_2^2 {d \choose m}$, whence, the first sum on the right hand side of \eqref{db3} is bounded above by
		\[        \sum_{0 \leq u < t/2} {t \choose u}^2 d^u u! ((t-u)!)^2\n{f}_2^2 {d \choose m}    =  E_0 (t!)^2 {d \choose m}  \sum_{0 \leq u < t/2} \frac{d^u}{u!} . \]
		For any $0 \leq u < t/2 = m \leq d/2$, observe that 
		\[ \frac{d^u}{u!} \leq \frac{1}{2^{t/2 - u}}\frac{d^{t/2}}{(t/2)!}. \]
		Thus the right hand side of the preceding inequality is bounded above by
		\[ E_0 (t!)^2 \frac{d^m}{m!} \frac{d^{t/2}}{(t/2)!} \sum_{0 \leq u < t/2} \frac{1}{2^{t/2 - u}} < E_0 t! {t \choose t/2} d^t \leq E_0 t! (2d)^t. \]
		Putting this together with \eqref{db3} and \eqref{db4} gives us
		\[ E_t \leq E_0 t! (2d)^t + 2dt E_{t-1}.\qedhere\]
	\end{proof}
	
	\begin{remark}\label{rmk:hypercube}
		This proof of \cref{asymenergy} we have presented can be modified to prove  \cref{hammingcor}.\footnote{In fact, the argument we present here can be extended to prove a weighted version of \cref{hammingcor} where $T(X)$ is replaced by the sum of $f(x)f(y)$ over pairs $x,y\in X$ at Hamming distance $2$ and the bound $\abs{X}\leq\binom dm$ is replaced by $\n{f}_1^2/\n{f}_2^2\leq\binom dm$. For notational simplicity, we only treat the case $f=\1_X$ in this remark.} 
		Given $y\in(\Z/2\Z)^d$, write
		\[\mathcal E_s(y):=\{(i_1,\ldots,i_{2s})\in[d]^{2s}:y=e_{i_1}'+\cdots+e_{i_{2s}}'\},\]
		and define $\mathcal R_s(y)$ to be the set of elements $\tau:=(i_1,\ldots,i_{2s})\in\mathcal E_s(y)$ with $i_1,\ldots,i_{2s}$ distinct. 
		A \emph{cancelling pair} in $\tau$ becomes any pair of equal entries. Defining $E_s$ analogously to \eqref{db2}, the inequality \cref{hldr} can be proven in precisely the same fashion, and we have
		\begin{equation}\label{eq:E1T}
			E_1=\sum_{x_1,x_2\in X}\abs{\mathcal E_1(x_1-x_2)}=2T(X)+d\abs{X}.
		\end{equation}
		Defining $\mathcal E_s^u(y)$ similarly, the analogue of \cref{claim57} in this setting reads
		\begin{equation}\label{eq:lowering-hyp}
			\abs{\mathcal E_s^u(y)}\leq ds(2s-1)\abs{\mathcal E_{s-1}^{u-1}(y)}/u.
		\end{equation}
		The extra factor of $2-1/s$ in \eqref{eq:lowering-hyp} as compared to \cref{claim57} is due to the fact that there are $\binom{2s}2=s(2s-1)$ places in which to insert a cancelling pair into an element of $\mathcal E_{s-1}^{u-1}(y)$, rather than $s^2$. Similarly to \cref{claim4}, one can use \eqref{eq:lowering-hyp} to prove
		\[E_t\leq 2(2t-1)dE_{t-1}+2\binom dm\frac{d^{t/2}(2t)!}{2^{t/2}(t/2)!}E_0\]
		for $t=2m$, whenever $m\leq d/4$.
		When $m>1$, the above inequality is enough to obtain the bound $E_1\leq 8mdE_0$. When $m=1$, one can prove such a bound by hand using \eqref{eq:E1T} and $T(X)\leq\abs{X}^2$, while for $m>d/4$ such a bound follows trivially from \eqref{trivialbound}.
	\end{remark}

\section{Theorem \ref{asymenergy} and Hypercontractivity}\label{sec:hypercontractivity}

In this section, we will consider the following corollary of Theorem \ref{asymenergy}.

\begin{prop}\label{prop:weighted}
	Let $d\in\N$, let $S=\{e_1,\ldots,e_d\}$ be the canonical basis of $\Z^d$, and let $X\subset\Z^d$ be a non-empty set such that $\abs{X}\leq\binom dm$ for some $m\in\N$. Then, for any $f\colon X\to[0,\infty)$, we have
	\[E_f(X,S)\leq 16 md\cdot\n{f}_2^2.\]
\end{prop}

This follows from Theorem \ref{asymenergy} in a straightforward manner. 

\begin{proof}[Proof of Proposition \ref{prop:weighted}]
	Cauchy's inequality implies that any function $f \colon \mathbb{R}^d \to [0, \infty)$ which is supported on some finite, non-empty set $X \subseteq \mathbb{R}^d$ also satisfies
	\[  \n{f}_1^2 \leq |X| \n{f}_2^2 .\]
	Thus if $|X| \leq {d \choose m}$, the hypothesis \eqref{ell12} of Theorem \ref{asymenergy} is automatically satisfied. Applying Theorem \ref{asymenergy} immediately gives us the desired conclusion.
\end{proof}

Our main aim of this section is to give an alternative proof of Proposition \ref{prop:weighted} which utilises ideas from the vast collection of results associated with the \emph{Hypercontractivity} phenomenon. We begin by introducing some relevant notation. Let $G = (\Z/2\Z)^d$. Given any function $f\colon G \to \R$ and any $1 \leq p < \infty$, define\footnote{ To fit with the standard conventions of the hypercontractivity literature, we use scaling conventions here which differ from those in the rest of the paper.}
\[  \n{f}_{L^p} = \left( 2^{-d} \sum_{x \in G } |f(x)|^p \right)^{1/p} \ \ \text{and} \ \ \n{f}_{L^\infty} = \max_{x \in G} |f(x)|.  \]
For any two functions $f, g \colon G \to \R$, we define their convolution $f * g \colon G \to \R$ by setting
\[ (f *g )(x) = 2^{-d} \sum_{y \in G} f(y) g(x-y)  = 2^{-d} \sum_{y \in G} f(y) g(x+y)\]
for all $x \in G$. We set $\langle f,g \rangle = (f*g)(0)$.

Next, for any $0 \leq \rho \leq 1$, define the function $\mu_{\rho} \colon G \to [0 ,\infty)$ as
\[ \mu_{\rho}(x_1, \dots, x_d)  = \prod_{1 \leq i \leq d} ( (1 + \rho) \1_{x_i = 0} + (1 - \rho) \1_{x_i = 1} )\]
for all $(x_1, \dots, x_d) \in G$. Finally, define the \emph{noise operator} $T_{\rho}$ by setting $T_{\rho} f = f * \mu_{\rho}$ for all $f\colon G \to \R$. Thus for  any $x \in G$, one has
\[ (T_{\rho}f)(x) = 2^{-d} \sum_{y \in G} \mu_{\rho}(y) f(x+y). \]

\begin{example}
	Let $\rho = 0$. Then $\mu_0(x) = 1$ for all $x \in G$. Thus, for any $f\colon G \to \R$, we see that $T_0 f$ is the constant function satisfying $(T_0 f)({x}) = 2^{-d}\sum_{y \in G}f(y)$ for all $x\in G$. Moreover, when $\rho  =1$, then $\mu_{1}(x) = 2^d$ for $x = 0 \in G$ and $\mu_1(x) = 0$ otherwise. This means that $T_1 f = f$ for all $f\colon G \to \R$.
\end{example}

Noting the above example, it is natural to ask, as $\rho$ varies in $[0,1]$, how much ``smoother'' $T_{\rho}(f)$ is compared to $f$. A straightforward application of H\"{o}lder's inequality gives us \emph{contractivity}, that is, for any $0 \leq \rho \leq 1$ and any $1 \leq p \leq \infty$ and any $f \colon G \to \R$, we have
\begin{equation} \label{contractivity}
	\n{T_{\rho} f}_{L^p} \leq \n{f}_{L^p} .
\end{equation}
The \emph{hypercontractivity} phenomenon suggests that $T_{\rho}f$ exhibits much more smoothness. In particular,  for any $1 \leq p \leq q \leq \infty$ and any $0 \leq \rho \leq (\frac{p-1}{q-1})^{1/2}$ and any $f \colon G \to \R$, one has
\begin{equation} \label{hypercontractivity2}
	\n{T_{\rho} f}_{L^q} \leq \n{f}_{L^p} .
\end{equation}
This immediately implies \eqref{contractivity} by a simple application of H\"{o}lder's inequality. Inequality \eqref{hypercontractivity2} seems to have been first proven by Bonami \cite{Bo1970}, with later developments  given by Beckner \cite{Bec1975} and Gross \cite{Gr1975}; see \cite[Chapter 9]{Od2014} for further references. 

We will now highlight the connection between such results and upper bounds on additive energies. Recall the following notation from \S\ref{unweightedsection}; for every $1 \leq i \leq d$, the element $e_i'\in G$ denotes the vector which has $1$ in its $i^{\rm th}$ coordinate and $0$ in all other coordinates. In particular, $\{e_1', \dots, e_d'\}$ is a basis of $G$. Observe that for any $0 \leq \rho \leq 1$ and any $X \subseteq G$ and any $f \colon X \to [0, \infty)$, one has
\begin{align}  \label{con1}
	\langle T_{\rho} f, f\rangle  & = 2^{-d}  \sum_{x \in G} (T_{\rho} f)(x) f(x) = 2^{-2d}\sum_{x,y \in G} \mu_{\rho}(y) f(x+y) f(x) \nonumber \\
	& \geq 2^{-d-1}  \bigg(\frac{1 + \rho}{2}\bigg)^{d-2}\bigg(\frac{1 - \rho}{2}\bigg)^2 \sum_{x\in X}\sum_{1\leq i \neq j \leq d} f(x) f(x + e_i' - e_j') \nonumber \\
	& = 2^{-d-1}  \bigg(\frac{1 + \rho}{2}\bigg)^{d-2}\bigg(\frac{1 - \rho}{2}\bigg)^2\left( E_f(X,\{e_1', \dots, e_d'\}) - d \sum_{x \in X} f(x)^2\right).  
\end{align}
On the other hand, hypercontractivity-type results naturally deliver  upper bounds for $\langle T_{\rho} f, f\rangle$. For our purposes, the classical hypercontractivity inequality \eqref{hypercontractivity2} does not seem strong enough to prove  \cref{prop:weighted}. We will instead use sharper hypercontractivity-type results given by Kirshner--Samorodnitsky \cite{KS2021} to prove the following.

\begin{lemma}\label{lem:ks-gen}
	Let $f\colon G \to[0,\infty)$ be supported on a set $X$ satisfying $|X| \leq \binom dm$ for some integer $1 \leq m \leq d/2$. Then
	\[E_f(X,\{e_1',\ldots,e_d'\})\leq 16 md \sum_{x \in X} f(x)^2.\]
\end{lemma}
\noindent We remark that, for $f=\1_X$, we have $E_f(X,\{e_1',\ldots,e_d'\})=d\abs{X}+2T(X)$, where $T$ is as in \cref{hammingcor}.

\cref{lem:ks-gen} allows us to prove Proposition \ref{prop:weighted} in a straightforward manner.

\begin{proof}[Proof of \cref{prop:weighted}]
    For $m>d/2$ the result is trivial, so we assume $m\leq d/2$.
	Let $\pi$ be the canonical projection map from $\Z^d$  to $G=(\Z/2\Z)^d$ and   let $g\colon G \to[0,\infty)$ be given by setting
	\[g(y)=\Big(\sum_{x\in\pi^{-1}(y)}f(x)^2\Big)^{1/2}\]
	for all $y \in G$. Thus, we have
	\[ \sum_{x \in X} f(x)^2 = \sum_{ y\in \pi(X)} g(y)^2 .\]
	The function $g$ is supported on $\pi(X)$ which satisfies $|\pi(X)| \leq |X| \leq \binom dm$, so \cref{lem:ks-gen} gives us
	\begin{equation} \label{niceinequality}
		E_g(\pi(X),\{e_1',\ldots,e_d'\})\leq 16 md \sum_{y \in \pi(X)} g(y)^2 = 16 md \sum_{x \in X} f(x)^2.
	\end{equation}
	Observe that
	\begin{align*}
		E_f(X,S)- d\sum_{x \in X} f(x)^2 & = \sum_{i\neq j}\sum_{x\in X}f(x)f(x+e_i-e_j) \\
		& = \sum_{i\neq j}\sum_{y\in\pi(X)}\sum_{x\in X\cap\pi^{-1}(y)}f(x)f(x+e_i-e_j).
	\end{align*}
	Applying Cauchy's inequality, we see that
	\begin{align*}
		E_f(X,S) &\leq d\sum_{x \in X} f(x)^2+\sum_{i\neq j}\sum_{y\in\pi(X)}\Big(\sum_{x\in\pi^{-1}(y)}f(x)^2\Big)^{1/2}\Big(\sum_{x\in\pi^{-1}(y)}f(x+e_i-e_j)^2\Big)^{1/2}\\
		&=d\sum_{y \in\pi(X)} g(y)^2 +\sum_{i\neq j}\sum_{y\in\pi(X)}g(y)g(y+e_i'-e_j')=E_g(\pi(X),\{e_1',\ldots,e_d'\}).
	\end{align*}
	Combining this with \eqref{niceinequality}  gives us the desired estimate.
\end{proof}

Hence, we will now focus on proving Lemma \ref{lem:ks-gen}. In this endeavour, we will require some further definitions from \cite{KS2021}. Thus, let $H\colon(0,1)\to[0,\infty)$ be the binary entropy function
\[H(\sigma)=-\sigma\log_2\sigma-(1-\sigma)\log_2(1-\sigma).\]
Let $\phi\colon[0,1/2]^2\to\R$ be the function as defined at the beginning of \cite[Section~2.A.6]{KS2021}. \cite[Lemma~21]{KS2021} implies that
\begin{align}\label{eq:KS-phi}
	\min_{0<\eps\leq 1/2}\big(\phi(\sigma,\eps)+1-H(\sigma)  -y\log_2\eps-& (1-y) \log_2(1-\eps)\big)
	\nonumber \\
	& =\sigma H\left(\frac y{2\sigma}\right)+(1-\sigma)H\left(\frac y{2(1-\sigma)}\right)
\end{align}
for $0\leq\sigma\leq1/2$ and $0\leq y\leq 2\sigma(1-\sigma)$. (The identity \eqref{eq:KS-phi} is the only property of the function $\phi$ we will use.) With this in hand, we now record the following result of  Kirshner--Samorodnitsky \cite[Corollary~28]{KS2021}.

\begin{lemma} \label{lem:ks}
	Let $0<\rho<1$ and $0\leq\sigma\leq1/2$, and let $\eps=(1-\rho)/2$.
	For any $f\colon\{0,1\}^d\to\R$ supported on a set of size at most $2^{dH(\sigma)}$, we have\footnote{Note that our convention for the parameter $\rho$ in the noise operator $T_\rho$ differs from theirs.}
	\[\langle T_\rho f,f\rangle\leq2^{d(\phi(\sigma,\eps)+1-H(\sigma))}\n{f}_{L^2}^2.\]
\end{lemma}

\begin{remark} \label{fourierremark}
	It is perhaps worth mentioning that despite the fact that the above inequality pertains to two objects defined in the physical space, Kirshner--Samorodnitsky prove Lemma \ref{lem:ks} by applying a Parseval-type identity to pass from $\langle T_{\rho} f, f\rangle$ to an inner product defined over the Fourier space and applying sharper \emph{level-$k$}-type inequalities there, which themselves are a consequence of a sharper hypercontractivity-type theorem.
\end{remark}

We are now ready to prove Lemma \ref{lem:ks-gen}.

\begin{proof}[Proof of Lemma \ref{lem:ks-gen}]
    For $m=1$, we may simply compute
    \[E_f(X,\{e_1',\ldots,e_d'\})\leq d\sum_{x\in X}f(x)^2+2\sum_{x,y\in X}f(x)f(y)\leq 3d\sum_{x\in X}f(x)^2,\]
    which is a strong enough bound.
    Henceforth, assume $m>1$.
	Let $\sigma=m/d$ and note that $\abs X\leq\binom dm\leq 2^{dH(\sigma)}$.
	We thus conclude from \cref{lem:ks}, for any $0<\rho<1$, that
	\[\langle T_\rho f,f\rangle\leq 2^{d(\phi(\sigma,\eps)+1-H(\sigma))}\n{f}_{L^2}^2,\]
	where $\eps=(1-\rho)/2$. Combining this with \eqref{con1} gives us
	\[E_f(X,\{e_1',\ldots,e_d'\})\leq \left(d+ 2\cdot 2^{d(\phi(\sigma,\eps)+1-H(\sigma)-t\log_2\eps-(1-t)\log_2(1-\eps))}\right)\sum_{x \in X} f(x)^2\]
	where $t=2/d$.
	Using \eqref{eq:KS-phi} to choose $\eps$ appropriately (which we may do since $\sigma=m/d\leq1/2$ and $t=2/d\leq \sigma\leq 2\sigma(1-\sigma)$), we thus have
	\[E_f(X,\{e_1',\ldots,e_d'\})\leq\left(d+ 2^{dr+1}\right) \sum_{x \in X} f(x)^2,\]
	where
	\[r=\sigma H\left(\frac1{d\sigma}\right)+(1-\sigma)H\left(\frac1{d(1-\sigma)}\right)=\frac mdH\left(\frac1{m}\right)+\frac{d-m}d H\left(\frac1{d-m}\right).\]
	The result follows, as in \cite[Proof~of~Corollary~7]{KS2021}, from observing that 
	\[ tH(1/t)\leq \log_2t+1/(\log 2) \]
	for $t\geq 1$. Indeed, the latter implies that $2^{dr}\leq e^2m(d-m)$, and so
	\[E_f(X,\{e_1',\ldots,e_d'\})\leq\left(d+ 2e^2m(d-m)\right)\sum_{x \in X} f(x)^2 \leq 16md\sum_{x \in X} f(x)^2.\qedhere\]
\end{proof}

\bibliographystyle{amsbracket}
\providecommand{\bysame}{\leavevmode\hbox to3em{\hrulefill}\thinspace}

\appendix

\section{A geometric approach to Theorem \ref{thm:main}} \label{appendix1}

In this section, we present a significantly digested proof of the following explicit version of \cref{thm:main} which was found by GPT-6 Astra when the main body of the document was in the proofreading phase.
In what follows, we call an affine-independent set of $d+1$ vectors in $\R^d$ an \emph{affine basis}.

\begin{theorem}\label{thm:explicit}
	Given a finite set $A\subset\R^d$ with full affine span, there exists some affine basis $S\subset A$ of $\R^d$ for which
	\[E(A,S)\leq(2d+3)\abs{A}-(2d+2).\]
\end{theorem}

\noindent \cref{thm:explicit} is tight when $A$ is the simplex. In this case, the only affine basis contained in $A$ is $A$ itself, and we have $E(A,A)=(2d+1)(d+1)$. We do not know if the multiplicative constant $2$ (for large $d$ and $\abs{A}$) is the best possible.

Unlike the proof we presented in the main body of the document, the proof of \cref{thm:explicit} is geometric in nature.
It follows the following three steps:
\begin{enumerate}
	\item Find a norm $\n{\cdot}$ on $\R^d$ which is scaled appropriately to capture the ``orientation'' of the set $A$.\footnote{This norm will be the Euclidean norm with respect to some basis of $\R^d$.}
	
	\item Show that, if a vector $v$ is large with respect to the norm $\n{\cdot}$, then $A\cap (A+v)$ cannot be large.
	
	\item Choose a probability distribution on affine bases contained within $A$ which biases against choosing vectors which are nearby with respect to $\n{\cdot}$.
\end{enumerate}
Once these steps are accomplished, one may simply choose $S\subset A$ an affine basis from the distribution constructed in (3) and upper-bound the expected value of $E(A,S)$ using (2). 

We begin the proof by defining the notion which captures (1).
\begin{defn}
	A finite, nonempty set $A\subset\R^d$ is in \emph{isotropic position} if
	\[\frac1{\abs{A}}\sum_{a\in A}a=0\quad\text{and}\quad\frac1{\abs{A}}\sum_{a\in A}aa^\intercal=I_d.\]
\end{defn}
This notion is standard in the convex geometry literature; for a reference on the use of isotropic position, see the book \cite{isotropicbook}.

\begin{lemma}\label{lem:frame-select}
	Let $A\subset\R^d$ be finite and nonempty.
	Suppose that the affine span of $A$ is all of $\R^d$.
	Then there is an affine transformation $\rho$ such that the image $\rho(A)$ of $A$ is in isotropic position.
\end{lemma}

In other words, given a full-dimensional subset $A$ of an abstract $d$-dimensional Euclidean space, one can assign a zero vector and a coordinate frame which makes $A$ isotropic.

\begin{lemma}\label{lem:inter}
	Let $A$ be in isotropic position and let $v\in\R^d$ be non-zero. We have
	\[\abs{A\cap (A+v)}\leq\frac{2\abs{A}}{\n{v}^2}.\]
\end{lemma}

\begin{lemma}\label{lem:determinant-dist}
	Let $A$ be a finite subset of $\R^d$ in isotropic position.
	There exists a probability distribution $\mathcal P$ supported on affine-independent subsets of $A$ of size $d+1$ such that, for any distinct $x,y\in A$,
	\[\Pr_{S\sim\mathcal P}[x,y\in S]=\frac{(1+\n{x}^2)(1+\n{y}^2)-(1+(x\cdot y))^2}{\abs{A}^2}.\]
\end{lemma}

We will prove and discuss these lemmas in the following subsections. 
First, we deduce \cref{thm:explicit}.

\begin{proof}[Proof of \cref{thm:explicit}]
	By \cref{lem:frame-select}, and the fact that additive energy is invariant under affine transformation, we may assume $A$ is in isotropic position.
	Let $\mathcal P$ be the probability distribution on affine bases $S\subset A$ given by \cref{lem:determinant-dist}. We compute
	\begin{align*}
		\EE_{S\sim\mathcal P} E(A,S)
		&=\sum_{x,y\in A}\Pr[x,y\in S]\cdot\abs{(A+x)\cap(A+y)}\\
		&=(d+1)\abs{A}+\sum_{\substack{x,y\in A\\x\neq y}}\Pr[x,y\in S]\cdot\abs{A\cap (A+x-y)}\\
		&\leq(d+1)\abs{A}+\sum_{\substack{x,y\in A\\x\neq y}}\frac{(1+\n{x}^2)(1+\n{y}^2)-(1+(x\cdot y))^2}{\abs{A}^2}\cdot\frac{2\abs{A}}{\n{x-y}^2},
	\end{align*}
	where in the last line we have used both \cref{lem:inter} and the conclusion of \cref{lem:determinant-dist}.
	Observe that we have the bound
	\begin{align*}
		(1+\n{x}^2)(1+\n{y}^2)-(1+(x\cdot y))^2
		&=1+\n{x}^2+\n{y}^2+\n{x}^2\n{y}^2-1-2(x\cdot y)-(x\cdot y)^2\\
		&=\n{x-y}^2+\n{x}^2\n{y}^2-(x\cdot y)^2\\
		&\leq\n{x-y}^2+\left(\frac{\n{x}^2+\n{y}^2}2\right)^2-(x\cdot y)^2\\
		&=\n{x-y}^2\left(1+\n*{\frac{x+y}2}^2\right).
	\end{align*}
	Thus, we see that
	\begin{align*}
		\EE_{S\sim\mathcal P}E(A,S)
		&\leq(d+1)\abs{A}+\frac2{\abs{A}}\sum_{\substack{x,y\in A\\x\neq y}}\left(1+\n*{\frac{x+y}2}^2\right)\\
		&=(d+1)\abs{A}+2(\abs{A}-1)+\frac1{2\abs{A}}\sum_{x,y\in A}(\n{x}^2+2(x\cdot y)+\n{y}^2)-\frac2{\abs{A}}\sum_{x\in A}\n{x}^2.     \end{align*}
	Since $A$ is in isotropic position, we have that 
	\[ \sum_{x\in A}\n{x}^2=\tr\sum_{x\in A}xx^\intercal=\abs{A}d \ \ \text{and} \ \ \sum_{x,y\in A}x\cdot y=0 .\] 
	This combines with the preceding inequality to give us
	\[\EE_{S\sim\mathcal P}E(A,S)  \leq (d+3)\abs{A}-2+d\abs{A}-2d=(2d+3)\abs{A}-(2d+2). \]
	Taking any $S$ with $E(A,S)\leq\EE_{S\sim\mathcal P}E(A,S)$ gives the result.
\end{proof}

\begin{remark}
	The saving in this proof \cref{thm:main} over the trivial bound by a factor of $d$ comes primarily from \cref{lem:inter}. 
    Indeed, if $A$ is isotropic, one may compute
    \[\EE_{x,y\in A}\n{x-y}^2=2d,\]
    where the expectation is over uniformly random points in $A$. So, one might hope that, for a ``typical'' pair of points $x,y\in A$, we have $\n{x-y}^2\gg d$ and thus (by \cref{lem:inter}) $\abs{A\cap (A+y-x)}\ll\abs{A}/d$.
    If this bound were to hold for all pairs of elements of some $S\subset A$ of size on the order of $d$, we would have $E(A,S)\ll\abs{A}d$.
    
    \cref{lem:determinant-dist} provides a particularly clean way of finding (enough) such pairs simultaneously, but there are other approaches as well.
    One alternative proof (which was also found by GPT-6 Astra)
    applies the additional computation that, for any set $A \subset \R^d$ in isotropic position and any affine subspace $L \subseteq \mathbb{R}^d$ of dimension $r < d$, one has 
	\begin{equation}\label{eq:far-from-subspace}
    \frac1{\abs{A}} \sum_{a \in A} \left( \inf_{x \in L} \n{a - x}^2 \right) \geq d - r. 
    \end{equation}
    (That is, isotropic sets do not cluster too strongly around affine subspaces.)
	The bound \eqref{eq:far-from-subspace} can be used to greedily construct an affine-independent set $S \subseteq A$ of size $\lceil d/3 \rceil$ such that for any distinct $s_1, s_2 \in S$, one has $\n{s_1 - s_2}^2 > 2d/3$: at each step, simply add the point in $A$ which is furthest away from the affine span of the previously chosen points. Applying \cref{lem:inter} leads to the bound $E(A,S) < 3|A||S|$.
\end{remark}

\subsection{Attaining isotropic position}

In this subsection, we prove \cref{lem:frame-select}.
This lemma is certainly not new; see for example \cite[Proposition~2.3.3]{isotropicbook} in the convex setting.

\begin{proof}[Proof of \cref{lem:frame-select}]
	Let $m=\frac1{\abs{A}}\sum_{a\in A}a$ be the center of mass of $a$, and consider the matrix
	\[M=\frac1{\abs{A}}\sum_{a\in A}(a-m)(a-m)^\intercal.\]
	The matrix $M$ is symmetric and positive definite: if $v\in\R^d\setminus\{0\}$, then
	\[v^\intercal Mv=\frac1{\abs{A}}\sum_{a\in A}((a-m)\cdot v)^2.\]
	This quantity is positive, since if it were zero then every point $a\in A$ would lie in the hyperplane $\{x \in \R^d :(x-m)\cdot v=0\}$, contradicting the fact that $A$ is $d$-dimensional.     Since $M$ is positive definite, we can write $M=N^2$ for some positive definite matrix $N$. Then $B:=\{N^{-1}(a-m):a\in A\}$ satisfies
	\[\frac1{\abs{B}}\sum_{b\in B}b=N^{-1}\frac1{\abs{A}}\sum_{a\in A}(a-m)=0\]
	and
	\begin{align*}
		\frac1{\abs{B}}\sum_{b\in B}bb^\intercal
		&=\frac1{\abs{A}}\sum_{a\in A}N^{-1}(a-m)(a-m)^\intercal N^{-\intercal}\\
		&=N^{-1}\left(\frac1{\abs{A}}\sum_{a\in A}(a-m)(a-m)^\intercal\right)N^{-1}=N^{-1}MN^{-1}=I_d.
	\end{align*}
	So, the affine transformation $\rho(x)=N^{-1}(x-m)$ puts $A$ into isotropic position.
\end{proof}

\subsection{Intersecting with translates}

In this subsection, we prove \cref{lem:inter}. It will be useful to first isolate the following computation: if $A$ is in isotropic position, then
\begin{equation}\label{eq:dotsq}
	\sum_{a\in A}(a\cdot v)^2=v^\intercal\left(\sum_{a\in A}aa^\intercal\right)v=v^\intercal(\abs{A}I_d)v=\abs{A}\cdot\n{v}^2.
\end{equation}
Now, the idea behind the proof of \cref{lem:inter} is as follows: for a pair of points $(x,x+v)$ which both lie in $A$, the values of $x\cdot v$ and $(x+v)\cdot v$ differ by at least $\n{v}^2$, and thus contribute at least $(-\n{v}^2/2)^2+(\n{v}^2/2)^2=\n{v}^4/2$ to the sum \eqref{eq:dotsq}.
Thus, there cannot be too many such contributions.
Carrying this argument through naively loses a factor of two as compared to the statement of \cref{lem:inter}; by being slightly more careful about counting patterns of the form $\{x,x+v,x+2v,\ldots\}$ in $A$, we may prove \cref{lem:inter}.

\begin{proof}[Proof of \cref{lem:inter}]
	Consider the graph $G$ whose vertex set is $A$ and where elements $x$ and $x-v$, for any $x\in A\cap(A+v)$, are connected by an edge.
	The graph $G$ has $\abs{A\cap (A+v)}$ edges.
	The connected components of $G$ are paths of the form $x,x+v,\ldots,x+(m-1)v$ for various $x$ and $m$; let $X$ be the set of such pairs $(x,m)$. For $(x,m)\in X$, we have $x+mv\not\in A$. We conclude that
	\begin{equation}\label{eq:A-line-partition}
		A=\bigsqcup_{(x,m)\in X}\{x,x+v,\ldots,x+(m-1)v\};\qquad \abs{A\cap(A+v)}=e(G)=\sum_{(x,m)\in X}(m-1).
	\end{equation}
	Now, the subset of $A$ corresponding to $(x,m)\in X$ contributes to \eqref{eq:dotsq} depending on $m$:
	\begin{align*}
		\sum_{k=0}^{m-1}((x+kv)\cdot v)^2
		&=\sum_{k=0}^{m-1}(x\cdot v+k\n{v}_2^2)^2\\
		&\geq\n{v}^4\min_{t\in\R}\sum_{k=0}^{m-1}(t+k)^2\\
		&=\n{v}^4\sum_{k=0}^{m-1}\left(-\frac{m-1}2+k\right)^2=\n{v}^4\cdot\frac1{12}m(m^2-1).
	\end{align*}
	We conclude from \eqref{eq:dotsq} and \eqref{eq:A-line-partition} that
	\[\abs{A}\cdot\n{v}^2=\sum_{a\in A}(a\cdot v)^2=\sum_{(x,m)\in X}\sum_{k=0}^{m-1}((x+kv)\cdot v)^2\geq\n{v}^4\sum_{(x,m)\in X}\frac1{12}m(m^2-1).\]
	Using the inequality $\frac16m(m^2-1)\geq m-1$ for all $m \in \mathbb{N}$ and \eqref{eq:A-line-partition}, we obtain the bound
	\[\abs{A}\cdot\n{v}^2\geq\n{v}^4\sum_{(x,m)\in X}\frac12(m-1)=\frac12\abs{A\cap (A+v)}\cdot\n{v}^4.\qedhere\]
\end{proof}

\begin{remark}
	The multiplicative constant $2$ in \cref{lem:inter} is tight.
	For $d=1$ one can take $v=2e_1$ and $A$ to consist of $\abs{A}/2$ copies each of $e_1$ and $-e_1$ (perhaps translated very slightly so that $A$ is a set rather than a multiset); for larger $d$, simply add to $A$ a regular $d$-simplex in the hyperplane $e_1^\perp$ at the appropriate scale to make $A$ isotropic. 
\end{remark}

\begin{remark}
    The idea behind \cref{lem:inter} --- if $S$ is a set in isotropic position then the size of the intersection $S\cap(S+v)$ decays as $\n{v}_2$ grows --- is not new, although this specific statement seems to be.
    A similar statement in one dimension was proven by Nishiyama \cite[Theorem~1(a)]{Ni2022} in the context of total variation between probability distributions on $\R$: from this statement, one can extract the bound
    \[\abs{A\cap (A+v)}\leq\frac{4\abs{A}}{\n{v}_2^2+4}\]
    for $A$ in isotropic position, which is comparable to \cref{lem:inter} (and has the advantage of remaining nontrivial even for $v$ near zero).
    In higher dimensions, a version of such a statement also appeared in an early version of a preprint of the second author \cite{Sc2024}: here, a bound on $K\cap (K+v)$ decaying exponentially in $\n{v}_2$ was proven for a convex body $K$ in isotropic position. (Such a bound relies substantially on the convexity of $K$, but a bound decaying quadratically in $\n{v}_2$ is sufficient for our purposes.)
\end{remark}

\subsection{A distribution on affine bases}

Finally, we prove \cref{lem:determinant-dist}.
The probability distribution $\mathcal P$ constructed in the proof is a fairly classical example of a \emph{determinantal probability measure}, and all of the facts we need are contained, for example, in the survey article of Lyons \cite{Lyons} on the topic.
In particular, the distribution we use is discussed in \cite[Sections~2~and~5]{Lyons}.

Let $X$ be a finite set in $\R^n$.
Given any subset $Y\subset X$, write $M_Y$ for the $n\times\abs{Y}$ matrix whose columns are the elements of $Y$. 
Note that the matrix $M_Y$ is only defined up to permutation of the columns.
The associated matrix
\[M_YM_Y^\intercal=\sum_{x\in Y}xx^\intercal\in\R^{n\times n}\]
is independent of the ordering of the columns of $M_Y$, and is nonsingular whenever $Y$ spans $\R^n$.
Moreover, the Cauchy--Binet formula gives us
\begin{equation}\label{eq:cauchy-binet}
	\sum_{\substack{Y\subset X\\\abs{Y}=n}}\det(M_Y)^2=\sum_{\substack{Y\subset X\\\abs{Y}=n}}\det(M_YM_Y^\intercal)=\det(M_XM_X^\intercal).
\end{equation}
If $X$ spans $\R^n$ (so that $M_XM_X^\intercal$ is nonsingular), then the identity \eqref{eq:cauchy-binet} inspires the following probability distribution on subsets $Y\subset X$ of $X$ of size $n$: choose $Y\subset X$ with probability $\det(M_YM_Y^\intercal)/\det(M_XM_X^\intercal)$.
Call this distribution $\mathcal P_X$. 
The support of $\mathcal P_X$ is exactly the set of bases of $\R^n$ contained in $X$.
Moreover, in a suitable sense, $\mathcal P_X$ places higher weight on bases which are ``well spread-out'': if $Y$ consists of vectors which are close together, then $\det(M_Y)$, and consequently $\det(M_YM_Y^\intercal)$, will be small.

If $\mathcal P$ is a determinantal probability measure on subsets $Y$ of a ground set $X$, it turns out that one can explicitly compute the probability that a particular subset $T\subset X$ is contained within $Y\sim\mathcal P$.
In our setting of interest, this is essentially the content of \cite[Theorem~5.1]{Lyons}.

\begin{lemma}\label{lem:conditions}
    Let $X\subset\R^n$. Let $T\subset X$ be fixed, and sample $Y\sim\mathcal P_X$. 
	Then
	\[\Pr[T\subset Y]=\det\left(M_T^\intercal(M_XM_X^\intercal)^{-1}M_T\right).\]
\end{lemma}
\noindent Note that, if $\abs{T}=n$, this probability reduces to $\det(M_TM_T^\intercal)/\det(M_XM_X^\intercal)$, and \cref{lem:conditions} collapses to the definition of $\mathcal P_X$.

For completeness of exposition, we present a proof of \cref{lem:conditions} using the Cauchy--Binet formula and Schur complements, which is essentially equivalent to the proof given in \cite{Lyons}.

\begin{proof}[Proof of \cref{lem:conditions}]
    Let $k=\abs{T}$.
    Note that the definition of $\mathcal P_X$ commutes with linear transformations, that is, if $L\colon\R^n\to\R^n$ is an invertible linear map, the distribution found by sampling $Y\sim\mathcal P_X$ and applying $L$ to each element of $Y$ is the same as the distribution $\mathcal P_{LX}$.
    If $T$ is linearly dependent, then $\Pr[T\subset Y]=0$ while the matrix $M_T^\intercal(M_XM_X^\intercal)^{-1}M_T$ has rank strictly less than $\abs{T}$, and so there is nothing to prove.
    As such, we may assume $T$ is linearly independent (in particular, $k\leq n$), and thus a change of basis allows us to assume that $T=\{e_1,\ldots,e_k\}$.
    
    For each element $x\in X$, write $u_x\in\R^k$ (resp.\ $v_x\in\R^{n-k}$) for the first $k$ (resp.\ last $n-k$) entries of $x$.
    Write $Z:=\{v_x:x\in X\setminus T\}$.
    Consider some set $T\subset Y\subset X$ with $\abs{Y}=n$, and write $Y=T\sqcup T'$.
    The matrix $Y$ has $k$ columns $e_1,\ldots,e_k$ and $n-k$ columns of the form $u_x\oplus v_x$ for each $x\in T'$.
    As a result,
    \[\det(M_Y)^2=\det\begin{pmatrix}I_k&*\\ 0 &M_{\{v_x:x\in T'\}}\end{pmatrix}^2=\det\left(M_{\{v_x:x\in T'\}}\right)^2.\]
    This enables us to compute
    \begin{align*}
    \det(M_XM_X^\intercal)\cdot\Pr[T\subset Y]
    &=\sum_{\substack{T\subset Y\subset X\\\abs{Y}=n}}\det(M_Y)^2\\
    &=\sum_{\substack{T'\subset X\setminus T\\\abs{T'}=n-k}}\det(M_{T\sqcup T'})^2
    &=\sum_{\substack{T''\subset Z\\\abs{T''}=n-k}}\det(M_{T''})^2=\det\left(\sum_{v\in Z}vv^\intercal\right),
    \end{align*}
    where we have used the Cauchy--Binet formula in the last equality.

    Now, consider the block decompositions
    \[M_XM_X^\intercal=\begin{pmatrix}N_{00}&N_{01}\\N_{10}&N_{11}\end{pmatrix}\quad\text{and}\quad (M_XM_X^\intercal)^{-1}=\begin{pmatrix}\tilde N_{00}&\tilde N_{01}\\\tilde N_{10}&\tilde N_{11}\end{pmatrix},\]
    where $N_{00},\tilde N_{00}\in\R^{k\times k}$ and $N_{11},\tilde N_{11}\in\R^{(n-k)\times(n-k)}$. Our previous computation shows that
    \begin{equation}\label{eq:probcomp}
        \Pr[T\subset Y]=\det(M_XM_X^\intercal)^{-1}\det\left(\sum_{v\in Z}vv^\intercal\right)=\frac{\det N_{11}}{\det(M_XM_X^\intercal)}.
    \end{equation}
    On the other hand, we have
    \begin{equation}\label{eq:rhscomp}
        M_T^\intercal(M_XM_X^\intercal)^{-1}M_T=\begin{pmatrix}I_k&0\end{pmatrix}\begin{pmatrix}\tilde N_{00}&\tilde N_{01}\\\tilde N_{10}&\tilde N_{11}\end{pmatrix}\begin{pmatrix}I_k\\0\end{pmatrix}=\tilde N_{00}.
    \end{equation}
    Now, we compute 
    \begin{equation}\label{eq:schur}        \begin{pmatrix}N_{00}&N_{01}\\N_{10}&N_{11}\end{pmatrix}\begin{pmatrix}I_k& 0\\-N_{11}^{-1}N_{10}&I_{n-k}\end{pmatrix}=\begin{pmatrix}N_{00}-N_{01}N_{11}^{-1}N_{10}&N_{01}\\ 0 &N_{11}\end{pmatrix}.
    \end{equation}
    Taking determinants in \eqref{eq:schur} gives $\det(M_XM_X^\intercal)=\det N_{11}\det(N_{00}-N_{01}N_{11}^{-1}N_{10})$, while inverting \eqref{eq:schur} we find $\tilde N_{00}=(N_{00}-N_{01}N_{11}^{-1}N_{10})^{-1}$.
    That is,
    \[\det(M_XM_X^\intercal)=\det N_{11}\det(N_{00}-N_{01}N_{11}^{-1}N_{10})=\frac{\det N_{11}}{\det \tilde N_{00}}.\]
    The result now follows from \eqref{eq:probcomp} and \eqref{eq:rhscomp}.
\end{proof}

We now begin to apply \cref{lem:conditions} to our setting.
If $M_XM_X^\intercal=\lambda I_n$ for some scalar $\lambda$ and $T=\{x,y\}$ is a two-element set, we have
\begin{equation}\label{eq:2pt-conditions}
\det(M_T^\intercal M_X^{-\intercal}M_X^{-1} M_T)=\begin{vmatrix}\lambda^{-1}x^\intercal x & \lambda^{-1}x^\intercal y\\ \lambda^{-1}y^\intercal x & \lambda^{-1}y^\intercal y\end{vmatrix}=\lambda^{-2}\left(\n{x}^2\n{y}^2-(x\cdot y)^2\right).
\end{equation}
With this computation and \cref{lem:conditions} in hand, we can now prove \cref{lem:determinant-dist}.

\begin{proof}[Proof of \cref{lem:determinant-dist}]
    Writing $X=\{(a\oplus 1):a\in A\}\subset\R^{d+1}\}$, we have
    \[\sum_{x\in X}xx^\intercal=\begin{pmatrix}\sum_{a\in A}aa^\intercal&\sum_{a\in A}a\\\sum_{a\in A}a^\intercal&\sum_{a\in A}1\end{pmatrix}=\abs{A}\cdot I_{d+1}.\]
    Consider the distribution $\mathcal P_X$ supported on bases of $\R^{d+1}$ contained in $X$ of size $d+1$. By the map $a\oplus 1\mapsto a$, we can form a corresponding distribution $\mathcal P$ on subsets of $A$ of size $d+1$; since $a_1,\ldots,a_{d+1}$ are affine-independent if and only if $a_1\oplus 1,\ldots a_{d+1}\oplus 1$ are linearly independent, the support of $\mathcal P$ is exactly the set of affine bases of $\R^d$ contained in $A$.
    So, applying \cref{lem:conditions} with \eqref{eq:2pt-conditions}, we have for $x,y\in A$ that
    \begin{align*}
    \Pr_{S\sim\mathcal P}[x,y\in S]
    &=\Pr_{Y\sim\mathcal P_X}[x\oplus 1,y\oplus1\in Y]\\
    &=\frac{\n{x\oplus 1}^2\n{y\oplus 1}^2-((x\oplus 1)\cdot(y\oplus 1))^2}{\abs{A}^2}\\
    &=\frac{(1+\n{x}^2)(1+\n{y}^2)-(1+(x\cdot y))^2}{\abs{A}^2}.\qedhere
    \end{align*}
\end{proof}

\end{document}